\documentclass[oneside]{amsart}
\usepackage{amssymb}

\usepackage{amscd}

\theoremstyle{plain}
\newtheorem{theorem}{Theorem}[section]
\newtheorem{corollary}[theorem]{Corollary}
\newtheorem{lemma}[theorem]{Lemma}

\theoremstyle{definition}

\newtheorem{definition}[theorem]{Definition}
\theoremstyle{remark}
\newtheorem*{acknowledgement}{Acknowledgement}

\newtheorem{remark}[theorem]{Remark}

\input tcilatex
\numberwithin{equation}{section}
\newlength{\customskipamount}
\newlength{\customkern}
\newlength{\customleftmargin}
\input cyracc.def
\font\tencyr=wncyr10
\def\cyr{\tencyr\cyracc}

\begin{document}
\title[Poincar\'{e}'s proof of the so-called Birkhoff-Witt theorem]{Poincar\'{e}'s proof of the so-called \linebreak Birkhoff-Witt theorem}
\subjclass{01A55, 01A60, 16S30, 17B35}
\keywords{Universal enveloping algebra of a Lie algebra, the canonical map
(symmetrization), the Birkhoff-Witt theorem}
\author{}
\maketitle

\begin{center}
\begin{equation*}
\begin{tabular}{ccc}
Tuong Ton-That &  & Thai-Duong Tran \\ 
Department of Mathematics &  & Department of Mathematics \\ 
University of Iowa &  & University of Iowa \\ 
Iowa City, IA 52242 U.S.A. &  & Iowa City, IA 52242 U.S.A. \\ 
E-mail: tonthat@math.uiowa.edu &  & E-mail: ttran@math.uiowa.edu
\end{tabular}
\end{equation*}
\bigskip

\emph{In honor of the 100th birthday of the article, ``Sur les groupes
continus''}\bigskip
\end{center}

%TCIMACRO{
%\TeXButton{TeX field}{\begin{quotation}%
%}}%
%BeginExpansion
\begin{quotation}%
%
%EndExpansion
\emph{Abstract:} A methodical analysis of the research related to the
article, ``Sur les groupes continus\textit{''}, of Henri Poincar\'{e}
reveals many historical misconceptions and inaccuracies regarding his
contribution to Lie theory. A thorough reading of this article confirms the
precedence of his discovery of many important concepts, especially that of
the \textit{universal enveloping algebra} of a Lie algebra over the real or
complex field, and the \textit{canonical map} (\textit{symmetrization}) of
the symmetric algebra onto the universal enveloping algebra. The essential
part of this article consists of a detailed discussion of his rigorous,
complete, and enlightening proof of the so-called Birkhoff-Witt theorem.%
%TCIMACRO{
%\TeXButton{TeX field}{\end{quotation}%
%}}%
%BeginExpansion
\end{quotation}%
%
%EndExpansion
\bigskip

TITRE: La d\'{e}monstration de Poincar\'{e} du th\'{e}or\`{e}me de
Birkhoff-Witt.

R\'{E}SUM\'{E}: Une analyse m\'{e}thodique des travaux faits en connexion
avec l'article, ``Sur les groupes continus''\textit{, }de Henri Poincar\'{e}
r\'{e}v\`{e}le des erreurs historiques et des jugements injustes en ce qui
concerne sa contribution \`{a} la th\'{e}orie de Lie. Une \'{e}tude
approfondie de cet article confirme la pr\'{e}c\'{e}dence de sa
d\'{e}couverte de plusieurs concepts importants; notamment de \textit{%
l'alg\`{e}bre enveloppante universelle} d'une alg\`{e}bre de Lie sur le
corps r\'{e}el ou le corps complexe, et de \textit{l'application canonique} (%
\textit{la sym\'{e}trisation}) de l'alg\`{e}bre sym\'{e}trique sur
l'alg\`{e}bre enveloppante universelle. L'essentiel de cet article consiste
en un examen approfondi de sa d\'{e}monstration rigoureuse et compl\`{e}te
du th\'{e}or\`{e}me de Birkhoff-Witt.

\textit{Liste de quelques mots-cl\'{e}s}: L'alg\`{e}bre enveloppante
universelle, l'application canonique (sym\'{e}trisation), le
th\'{e}or\`{e}me de Birkhoff-Witt.

\section{Introduction\label{intro}}

In our research on the universal enveloping algebras of certain
infinite-dimensional Lie algebras we were led to study in detail the
original proofs of the so-called Birkhoff-Witt theorem (more recently,
Poincar\'{e}-Birkhoff-Witt theorem). This, in turn, led us to the
investigation of Poincar\'{e}'s contribution to Lie theory (i.e., the theory
of Lie groups, Lie algebras, and their representations). To our great
surprise we discovered many historical misconceptions and inaccuracies, even
in some of the classics written by the leading authorities on the subject.
This discovery has puzzled us for some time, and we have sought the opinions
of several experts in the field. Their answers together with our thorough
reading of several original articles on the subject shed some light on this
mystery. We were astounded to find out that Poincar\'{e} was given neither
credit for his fundamental discovery of the \emph{universal enveloping
algebra} of a Lie algebra over a field of characteristic zero, nor for\ his
introduction of the \textit{symmetrization map},\textit{\ }and only a
cursory and belated\ acknowledgment of his contribution to the so-called 
\textit{Birkhoff-Witt theorem}, of which he gave a rigorous, complete,
beautiful, and enlightening proof. Indeed, in two of the most exhaustive
treatises on universal enveloping algebra \cite{Co} and \cite{Di}
Poincar\'{e}'s work \cite{Po1} was not mentioned. In many authoritative
textbooks treating Lie theory such as \cite{Ch}, \cite{C-E}, \linebreak \cite
{Ku}, \cite{Ja}, \cite{Va}, \cite{Hu}, \linebreak \cite{Kn},$\ldots $,
Poincar\'{e}'s discovery of the universal enveloping algebra and the
symmetrization map was ignored. In some books his name was left off the
Birkhoff-Witt theorem, and his fundamental article \cite{Po1} was not even
quoted. In \textit{Encyclopaedia of Mathematics} \cite{En} under the rubric
``Birkhoff-Witt theorem'' it was written ``$\ldots $The first variant of
this theorem was obtained by H. Poincar\'{e}; the theorem was subsequently
completely demonstrated by E. Witt [1937] and G.D. Birkhoff\footnote{%
Actually Garrett Birkhoff (1911--1996), not G.D. (Birkhoff) which are the
initials of George David Birkhoff (1884--1944), the father of Garrett. This
inaccuracy only occurs in the translation, not in the original (Russian)
version of the Encyclopaedia. We are grateful to Professor Sergei Silvestrov
for elucidating this fact to us.} [1937]$\ldots $''. Clearly the author,
T.S. Fofanova, did not read carefully \cite{Po1}; otherwise she would have
realized that Poincar\'{e} had discovered and completely demonstrated this
theorem at least thirty-seven years before Witt and Birkhoff. Why such
iniquities can happen to one of the greatest mathematicians of all times who
published these results \cite{Po1} in one of the most prestigious scientific
journals, \textit{Trans.\ Camb.\ Philos.\ Soc.}, on the occasion of the
jubilee of another great mathematician, Sir George Gabriel Stokes, is a most
interesting mystery that we shall attempt to elucidate in this article. But
before beginning our investigation we want to make it clear that our
intention is to study thoroughly one of the most fundamental discoveries by
one of the greatest minds in order to understand how important ideas are
created, and not to rectify such iniquities, for such a task is doomed to
fail as the force of habit always prevails; a fact very perspicuously
expressed in the following excerpt from \cite[p.~186]{Gi}, ``$\ldots $%
l'Hospital's rule, Maclaurin's series, Cramer's rule, Rolle's theorem, and
Taylor's series are familiar terms to calculus students. Actually, only one
of these five mathematicians was the original discoverer of the result
attributed to him, and that man was Rolle. The person who popularizes a
result generally has his name attached to it, although later it may be
learned that someone else had originally discovered the same result. For
practical purposes names are not changed, but even so, the mistakes seem to
compensate for one another. Although Maclaurin was credited with a series he
did not discover, a rule which he did originate is now known as Cramer's rule%
$\ldots $''. Besides, Poincar\'{e} is a member of the elite group of
mathematicians to whom many important mathematical discoveries are
attributed; indeed, in \textit{Encyclopaedia of Mathematics} \cite{En} 18
rubrics are listed under his name. Curiously, under the heading
``Poincar\'{e} last theorem'' the editorial comments state that ``[this
theorem] is also known as the Poincar\'{e}-Birkhoff fixed-point theorem,''
and the author, M.I. Vo\u{\i}tsekhovski\u{\i}, wrote ``$\ldots $it was
proved by him in a series of particular cases but he did not, however,
obtain a general proof of this theorem.''\footnote{%
Vo\u{\i}tsekhovski\u{\i} continues, ``The paper was sent by Poincar\'{e} to
an Italian journal two weeks before his death, and the author expressed his
conviction, in an accompanying letter to the editor, of the validity of the
theorem in the general case.'' Indeed, on December 9, 1911, having some
presentiments that he might not live long, Poincar\'{e} wrote a moving
letter to Guccia, director and founder of the journal \emph{Rendiconti del} 
\emph{Circolo Matematico di Palermo} (cf. \cite[t.~2, p.~LXVII]{Po2}),
asking his opinion regarding what has become known as ``Poincar\'{e}'s Last
Geometric Theorem'' (see \cite[\S7.4.2, 169--174]{B-G}, for an English
translation of the letter and an excellent discussion of the theorem). Mr.\
Guccia readily accepted the memoir for publication and it appeared on March
10, 1912, just a few months before Poincar\'{e}'s death on July 17, 1912
(Sur un th\'{e}or\`{e}me de G\'{e}om\'{e}trie, \textit{Rendiconti del
Circolo Matematico di Palermo},\textit{\ }33\textbf{, }375--407 = Oeuvres
VII, 499--538). Ultimately it was G.D. Birkhoff who gave a complete proof of
this theorem (Proof of Poincar\'{e}'s geometric theorem, \textit{Trans.\
Amer.\ Math.\ Soc.}, 14 (1913)\textbf{,} pp.\ 14--22 = \textit{Collected
Mathematical Papers I}, 673--681) and of its generalization to $n$
dimensions (Une g\'{e}n\'{e}ralisation \`{a} $n$ dimensions du dernier
th\'{e}or\`{e}me de g\'{e}om\'{e}trie de Poincar\'{e}, \textit{C.R.\ Acad.\
Sci. Paris, }192\textbf{\ (}1931), pp. 196--198 = \textit{Collected
Mathematical Papers II}, 395--397).} Misnaming mistakes seem to compensate
one another after all.

In his book \cite{Be}, E.T. Bell, who called Poincar\'{e} ``the Last
Universalist'', considered the Last Geometric Theorem as Poincar\'{e}'s
``unfinished symphony'' and wrote ``$\ldots $And it may be noted that
Poincar\'{e} turned his universality to magnificent use in disclosing
hitherto unsuspected connections between distant branches of mathematics,
for example, between continuous groups and linear algebra''. This is exactly
the impression we had when reading his article, ``Sur les groupes continus''.

\section{Poincar\'{e}'s work on Lie groups\label{Poincare's work}}

To assess Poincar\'{e}'s contribution to Lie theory in general we use two
main sources \cite{Po2} and \cite{Po3} and investigate in depth the
references cited therein. We start with the article, ``Analyse des travaux
scientifiques de Henri Poincar\'{e}, faite par lui-m\^{e}me\textit{''}%
\footnote{%
In \cite{Po2} this article is listed as published by \textit{Acta Math}., 30
(1913), pp. 90--92. In fact, it never existed as such; the editors of
Poincar\'{e}'s collected work probably found the manuscript of the article
among his papers with his annotations regarding the journal and the date of
publication but due to World War I it appeared eventually in \cite{Po3}.
This remark extends to all discrepancies between the intended and actual
dates of publication of many of Poincar\'{e}'s works in \cite{Po3}, for
example, \emph{Rapport sur les travaux de M. Cartan}.}\textit{\ }which was
written by Poincar\'{e} himself in 1901 (\cite[3--135]{Po3}) at the request
of G. Mittag-Leffler (cf. ``Au lecteur\textit{''} \cite[1--2]{Po3}). It is
part of vol.\ 38\textbf{\ }of the journal, \textit{Acta Math}.\textit{,}
published in 1921 in memory of Henri Poincar\'{e}. (Actually, most of vol. 39%
\textbf{\ }published in 1923 is also devoted to Poincar\'{e}'s work). In the
third part of the above-mentioned article, Section XII (Alg\`{e}bre) and
Section XIII (Groupes Continus) are devoted to his contribution to Lie
theory. Actually, we think that under Poincar\'{e}'s impetus
finite-dimensional continuous groups are eventually called Lie groups.
Indeed, Poincar\'{e} expressed repeatedly his great admiration for Lie's
work in \cite{Po4} and \cite{Po1} and wrote in \textit{Rapport sur les
travaux de M. Cartan} (\cite[137--145]{Po3}, with a curious footnote: 
\textit{Acta Math}., 38, printed on August 11, 1914): ``Je commencerai par 
\emph{les groupes continus et finis}, qui ont \'{e}t\'{e} introduits par Lie
dans la science; le savant norv\'{e}gien a fait conna\^{i}tre les principes
fondamentaux de la th\'{e}orie, et il a montr\'{e} en particulier que la
structure de ces groupes depend d'un certain nombre de constantes qu'il
d\'{e}signe par la lettre $c$ affect\'{e}e d'un triple indice et entre
lesquelles il doit y avoir certaines relations$\ldots $ une des plus
importantes applications \emph{des groupes de Lie}$\ldots $''. So far as we
know this is the first time that the name \emph{Lie groups} was explicitly
mentioned.

Poincar\'{e}'s first encounter with Lie theory probably dated back to his
article \cite{Po5} and its generalization \cite{Po6}. The problem he
considered there can be phrased in modern language as follows:

For $X=(x_{1},\ldots ,x_{n})\in \mathbb{C}^{n}$ let $\mathrm{GL}_{n}(\mathbb{%
C})$, the general linear group of all $n\times n$ invertible complex
matrices, act on $\mathbb{C}^{n}$ via $(X,g)\rightarrow Xg$, $g\in \mathrm{GL%
}_{n}(\mathbb{C)}$. Let $F(X)$ denote a homogeneous form of degree $m$
(i.e., a homogeneous polynomial of degree $m$ in $n$ variables $%
(x_{1},\ldots ,x_{n})$), find the subgroup $G$ of $\mathrm{GL}_{n}(\mathbb{C}%
)$ which preserves the form $F$; i.e., $F(Xg)=F(X),\forall g\in G$.
Conversely, given a subgroup $G$ of $\mathrm{GL}_{n}(\mathbb{C})$ find all
homogeneous forms that are left invariant by $G$. This is precisely the
problem of polynomial invariants (cf.\ \cite{We}).

In \cite{Po5} and \cite{Po6} he found all cubic ternary (of three variables)
and quaternary (of four variables) forms that are preserved by certain
Abelian groups (which he called ``faisceau de substitutions''), and he also
extended this result to the non-Abelian case. Conversely, he exhibited
explicit groups that preserve quadratic and cubic ternary and quaternary
forms. For example, in \cite[239--241]{Po5} he found the subgroup of the
unipotent group 
\begin{equation*}
\left\{ \left( 
\begin{array}{lll}
1 & \alpha & \beta \\ 
0 & 1 & \gamma \\ 
0 & 0 & 1
\end{array}
\right) ;\alpha ,\beta ,\gamma \in \mathbb{C}\right\}
\end{equation*}
which preserves the quadratic form 
\begin{multline*}
\left[ 
\begin{array}{lll}
x_{1} & x_{2} & x_{3}
\end{array}
\right] \left[ 
\begin{array}{lll}
A_{1} & B_{3} & B_{2} \\ 
B_{3} & A_{2} & B_{1} \\ 
B_{2} & B_{1} & A_{3}
\end{array}
\right] \left[ 
\begin{array}{l}
x_{1} \\ 
x_{2} \\ 
x_{3}
\end{array}
\right] \\
=A_{1}x_{1}^{2}+A_{2}x_{2}^{2}+A_{3}x_{3}^{2}+2B_{1}x_{2}x_{3}+2B_{2}x_{1}x_{3}+2B_{3}x_{1}x_{2},\quad A_{i},B_{i}\in 
\mathbb{C};1\leq i\leq 3.
\end{multline*}
Conversely, he showed that all quadratic forms which are left invariant by
the full unipotent group defined above must satisfy a certain partial
differential equation. Lie theory also plays an important role in
Poincar\'{e}'s work on the conformal representation of functions of two
variables which leads to the theory of relativity. In \cite{Po7} and \cite
{Po8} he studied the group of linear transformations which leave the
Minkowski's metric $x^{2}+y^{2}+z^{2}-t^{2}$ invariant which he called the
homogeneous Lorentz group, or as H.A. Lorentz wrote in \cite{Lo} ``groupe de
relativit\'{e}'', and discovered the Poincar\'{e} group which is the
semidirect product of the four-dimensional translation group with the
homogeneous Lorentz group.

In his analysis of his scientific accomplishments Poincar\'{e} classified
his work in seven topics which range from ``Differential Equations'' to
``Philosophy of Science''. His accomplishments in any single one of these
areas would already make him famous. Indeed Sir George H. Darwin
(1845--1912), a physicist and son of the famous Charles Darwin (1809--1882),
wrote in 1909 that Poincar\'{e}'s celestial mechanics would be a vast mine
for researchers for half a century (\cite[p.~652]{By}). Under rubric number
three \textit{``}Questions diverses de Math\'{e}matiques pures\textit{'',}
Algebra, Arithmetic, Group theory, and Analysis Situs (combinatorial
topology) are listed together, with Lie groups as a subsection of Group
theory. This gives a false impression that he had only a slight interest in
the subject. In fact with the exception of his first article on continuous
groups \cite{Po4} the other three articles are quite long: \cite{Po1} (35
pages), \cite{Po9} (47 pages) and \cite{Po10} (60 pages). In all these
articles he not only conveyed to the reader his keen interest in the subject
but also some of the difficulties that preoccupied him over a ten-year
period.

It was Lie's third theorem that motivated Poincar\'{e} to write \cite{Po4}.
This theorem can be stated in \textit{his notations}\footnote{%
Poincar\'{e} used parentheses instead of brackets for the commutator
products. To avoid confusion we replace the parentheses with the more
conventional brackets. He also used the notation $c_{iks}$ instead of the
more convenient notation $c_{ik}^{s}$ for the structure constants. We do not
however replace this notation, which does not cause any confusion, to
preserve as much as possible Poincar\'{e}'s style and terminology.} as
follows:

If $\{X_{1},\ldots ,X_{r}\}$ is a system of infinitesimal transformations
(i.e., vector fields) which satisfy the equation 
\begin{equation}
\lbrack X_{i},X_{k}]=\sum_{s=1,\ldots ,r}c_{iks}X_{s}  \label{eq2.1}
\end{equation}
then the structure constants $c_{iks}$ must satisfy the relations 
\begin{equation}
c_{kis}=-c_{iks},  \label{eq2.2}
\end{equation}
\begin{equation}
\sum_{k=1}^{r}\left( c_{iks}c_{jlk}+c_{lks}c_{ijk}+c_{jks}c_{lik}\right)
=0\qquad (1\leq i,j,l,s\leq r),  \label{eq2.3}
\end{equation}
which follow immediately from the fact that the bracket $\left[ \,\cdot
\,,\,\cdot \,\right] $ is skew symmetric, and the Jacobi identity 
\begin{equation*}
\lbrack \lbrack
X_{a},X_{b}],X_{c}]+[[X_{b},X_{c}],X_{a}]+[[X_{c},X_{a}],X_{b}]=0.
\end{equation*}
Conversely, if the coefficients $c_{iks}$ satisfy Eqs.\ (\ref{eq2.2}) and (%
\ref{eq2.3}) then there exists a system of infinitesimal transformations
verifying Eq.\ (\ref{eq2.1}), and hence a group of transformations with $r$
parameters.

In \cite{Po1} he gave a different proof of this theorem, especially for the
case of Lie algebras with non-trivial centers. His approach consists of
reformulating Lie's construction of the \textit{adjoint group}, deriving the
differential equations associated with the\textit{\ parametric group}, and
then showing that these equations can be integrated. More specifically, if $%
X_{1},\ldots ,X_{r}$ form a basis of a Lie algebra of infinitesimal
transformations, and 
\begin{equation*}
V=\sum v_{i}X_{i},\qquad T=\sum t_{i}X_{i},
\end{equation*}
then one has the adjoint representation 
\begin{equation*}
T\rightarrow T^{\prime }=e^{-V}Te^{V},
\end{equation*}
where $T^{\prime }=\sum t_{i}^{\prime }X_{i}$ and $e^{tX}=\sum_{n=0}^{\infty
}\frac{t^{n}}{n!}X^{n}$. The image of the adjoint representation is then
called the adjoint group. By setting 
\begin{equation*}
e^{V}e^{T}=e^{W},\text{\qquad where }W=\sum w_{i}X_{i},
\end{equation*}
it follows that the $w$ are functions of the $v$ and $t$, or, in other words
the transformation $e^{T}$ transforms $e^{V}$ into $e^{W}$, i.e., the $v$
into the $w$; and the group thus defined in $r$ variables is called the
parametric group associated with the system $\{X_{1},\ldots ,X_{r}\}$. In 
\cite{Po9} and \cite{Po10} he studied in great detail these differential
equations and the isomorphism between the adjoint and parametric groups. His
research into expressions of $W$ as a function of $U$ and $V$ in the formula 
$e^{U}e^{V}=e^{W}$ resulted in a precursory form of the
Baker-Campbell-Hausdorff formula (see \cite{W} for a discussion regarding
Poincar\'{e}'s contribution to this theorem and, e.g., \cite[p.~114]{Va} or 
\cite[Chap.~II]{Bo1}, for a more modern proof of this theorem).

\section{Poincar\'{e}'s discovery of the universal enveloping algebra and
the so-called Birkhoff-Witt theorem\label{Poincare discovery}}

In this section we shall expound the main theme of this article, namely,
Poincar\'{e}'s precedence in the discovery of the universal enveloping
algebra and the so-called Birkhoff-Witt theorem. For this purpose we shall
examine in detail his article \cite{Po1}. As a general rule we try to adhere
faithfully to his exposition, notation, and style as much as possible. But
in order to make our point we shall insert some comments, prove some claims
which Poincar\'{e} considered self-evident but did not seem to be so obvious
to us, and integrate his work in the more modern framework of Lie theory. At
first reading \cite{Po1} seems to be hastily written, repetitive, and
sometimes cryptic, and this might explain why not many people have read it;
especially for the readers to whom French is not their first language. But
by a careful analysis of \cite{Po1} one must conclude without a shade of
doubt that Poincar\'{e} had discovered the concept of the universal algebra
of a Lie algebra and gave a complete and rigorous proof of the so-called
Birkhoff-Witt theorem. As we shall see, his entire proof of this theorem,
with the exception of the claim that we will state as Theorem \ref
{Theorem3.3}, is quite rigorous and modern in language. For these reasons we
will translate the parts in \cite{Po1} that are relevant to our discussion
for the benefit of the readers who are not familiar with French. But before
going into the details we shall elaborate on why we consider his proof very
enlightening. For example, his introduction of \emph{the symmetrization map}
is quite natural by observing that the most elementary ``regular'' (or
``symmetric'') polynomials are the linear polynomials and their powers. And,
as it turned out, all symmetrized polynomials are linear combinations of
those. From the symmetrization comes out naturally the notion of \emph{%
equipollence} which, in turn, leads to the notion of \emph{equivalence}, and
ultimately to the definition of the \emph{universal enveloping algebra of a
Lie algebra}. In the proof of the fact stated as Theorem \ref{theorem3.8}
below, he introduced the notion of ``chains'' and cleverly showed that he
can add \emph{more} chains to paradoxically \emph{reduce} the number of
basic chains, and this enables him to proceed by induction (on the degree of
the regular polynomials). This ingenious idea foreshadows some techniques
used in the modern theory of \emph{word problem}.

For the remainder of this section, in order to capture Poincar\'{e}'s vivid
flow of ideas we shall use the present tense to present his exposition. For
convenience, we shall discuss the universal enveloping algebra first. Let $%
X_{1},\ldots ,X_{n}$ be $n$ elementary operators (Poincar\'{e} thinks of
these operators as vector fields but never really uses this fact here). Let $%
\mathcal{L}$ be the Lie algebra over a field of characteristic zero $\mathbb{%
K}$ (Poincar\'{e} thinks of $\mathbb{K}$ as $\mathbb{R}$ or $\mathbb{C}$ but
all concepts and proofs remain identical) generated by these $n$ elementary
operators which constitute a basis for $\mathcal{L}$. Let $\mathcal{A}$
denote the \textit{non-commutative }algebra of polynomials in $n$ variables $%
X_{1},\ldots ,X_{n}$ with coefficients in $\mathbb{K}$. Consider the set of
all elements of $\mathcal{A}$ of the form 
\begin{equation}
P(XY-YX-[X,Y])Q,  \label{eq3.1}
\end{equation}
where $P$ and $Q$ are arbitrary polynomials in $\mathcal{A}$, and where $%
[X,Y]$ denotes the bracket product of $\mathcal{L}$. Define an \textit{%
equivalence relation} $\sim $ in $\mathcal{A}$ by declaring that an element $%
A\in \mathcal{A}$ is equivalent to $0$ if $A$ is a linear combination of
elements of the form (\ref{eq3.1}) for some $P$ and $Q$ in $\mathcal{A}$,
and $A\sim A^{\prime },A^{\prime }\in \mathcal{A}$, if $A-A^{\prime }\sim 0$%
. Then the quotient algebra (or residue ring) thus defined is now called the 
\textit{universal enveloping algebra }of $\mathcal{L}$. In fact, this can be
rephrased in modern language as follows:

Let $\mathcal{T}$ denote the tensor algebra over the underlying vector space
of $\mathcal{L}$, then $\mathcal{T}$ is isomorphic to $\mathcal{A}$ (see,
for example, \cite[Prop.~10, p.~423]{La}). Let $J$ denote the two-sided
ideal of $\mathcal{T}$ generated by the tensors $X\otimes Y-Y\otimes X-[X,Y]$
where $X,Y\in \mathcal{L}$. Then the associative algebra $\mathcal{U}=%
\mathcal{T}/J$ is called the universal enveloping algebra of $\mathcal{L}$
(cf., e.g., \cite[p.~12]{Bo2} or \cite[p.~22]{Bo3}). Under the isomorphism
between $\mathcal{T}$ and $\mathcal{A}$ the ideal $J$ corresponds to the
two-sided ideal of $\mathcal{A}$ spanned by all elements of the form (\ref
{eq3.1}). Actually this is exactly Harish-Chandra's approach to the
universal enveloping algebra in \cite{H-C1}. Unaware that Poincar\'{e} had
defined this notion in \cite{Po1}, Harish-Chandra wrote the following
footnote: ``\textit{This algebra is the same as the one considered by
Birkhoff }[1937]\textit{\ and Witt }[1937],\textit{\ though their method of
construction is different''. }Indeed, as we mentioned earlier,
Harish-Chandra's definition of the universal enveloping algebra is identical
to Poincar\'{e}'s. He also wrote: ``\textit{In view of this 1-1
correspondence between representations of }$\mathcal{L}$\textit{\ and }$%
\mathcal{U}$\textit{\ it is appropriate to call }$\mathcal{U}$\textit{\ the
general enveloping algebra of }$\mathcal{L}$\textit{''}, and in his
fundamental paper \cite{H-C2} on the role of the universal enveloping
algebra of a semisimple Lie algebra in Lie theory, published two years
later, he replaced the word ``general'' with ``universal'', probably under
the influence of Birkhoff's work on universal algebras. Thus we can conclude
that it was Harish-Chandra who named this algebra discovered by Poincar\'{e}
``universal enveloping algebra''.

In his very influential book \cite{Ch}, C. Chevalley, one of the world's
leading experts in Lie theory and a founding member of Bourbaki, attributed
to Harish-Chandra the following theorem: \textit{``There exists a one-to-one
correspondence }(\textit{but not multiplicative}!)\textit{\ between elements
of }$\mathcal{U}$\textit{\ and those of the symmetric algebra of }$\mathcal{L%
}$\textit{; also if }$\mathcal{L}$\textit{\ is the Lie algebra of a Lie
group }$G$\textit{, then }$\mathcal{U}$\textit{\ is isomorphic to the
algebra of right }(\textit{or left})\textit{\ invariant differential
operators over the algebra of analytic functions on }$G$\textit{'' }(see
Chapter 5, vol.\ III of \cite{Ch}, especially \S 6).\textrm{\ }Note that
Chevalley and Harish-Chandra were colleagues at Columbia University during
this period. We shall establish below that Poincar\'{e} had already
discovered this correspondence (usually called the \textit{canonical
isomorphism }or\textit{\ }the\textit{\ symmetrization map})\textit{\ }and
had defined $\mathcal{U}$ as the algebra of right (or left) invariant
differential operators on $G$. It is also interesting to note that Chevalley
also gives a proof of the Birkhoff-Witt theorem \cite[vol.~III, Prop.~1,
p.~163]{Ch} without mentioning the work of Birkhoff and Witt. Anyhow, many
authors seem not to acknowledge Poincar\'{e}'s discovery of the fundamental
notion of universal enveloping algebra; for example, in the encyclopaedic
work \cite{Di} Poincar\'{e}'s work is not even referred to with regard to
this algebra.

Finally we are coming to the main part of this article, namely,
Poincar\'{e}'s proof of the Birkhoff-Witt theorem. But before going into
detail about the proof, we shall make some historical remarks. As discussed
earlier, none of the leading experts in Lie theory seemed to be aware of the
existence of Poincar\'{e}'s work on the universal enveloping algebra and his
proof of the so-called Birkhoff-Witt theorem prior to about 1956. Garrett
Birkhoff and Ernst Witt certainly didn't mention Poincar\'{e}'s work in \cite
{Bi} and \cite{Wi}, respectively, which both appeared in 1937. M. Lazard
generalized in \cite{Lz-1} and \cite{Lz-2} this theorem, which he called the
Witt theorem, but did mention \cite{Bi} and the work of Kourotchkine [1951].
So far as we know the authors who first noticed the proof of the
Birkhoff-Witt theorem already appeared in \cite{Po1} were H. Cartan and S.
Eilenberg, whose book \cite{C-E} appeared in 1956. Curiously, they called
the theorem the Poincar\'{e}-Witt theorem and did not refer to Birkhoff's
proof; moreover, they attributed the complete proof of the theorem to Witt.
It appears that Bourbaki was the first to call this theorem the
Poincar\'{e}-Birkhoff-Witt theorem in \cite{Bo3}, a recognition acknowledged
in arguably the most influential book on Lie algebras in the English
language \cite{Ja}. From then on this is the most prevalent name used for
this theorem; however, many authors of serious books on Lie theory such as 
\cite{Ku}, \cite{Co} and more recently \cite{Kn}, etc., continue to call it
the Birkhoff-Witt theorem.

Now let us carefully examine \cite{Po1}, especially the portion relevant to
our investigation, pp.\ 224--232.

The section heading is ``Calcul des polyn\^{o}mes symboliques''. Let $%
X,Y,Z,T,U,\ldots ,$ be $n$ \textit{elementary operators} (i.e., a basis for
a Lie algebra over a commutative field $\mathbb{K}$ of characteristic zero).
Consider the algebra of \textit{symbolic }(\textit{or formal})\textit{\
non-commutative polynomials} in these operators with coefficients in $%
\mathbb{K}$. Then as previously mentioned, we may identify this algebra with
the tensor algebra $\mathcal{T}$.

\begin{definition}
\label{Def3.1}Two monomials are said to be \textit{equipollent} if they
differ only by the order of their factors. This definition extends obviously
to two polynomials that are sums of pairwise equipollent monomials. Ex. $%
XY^{2}$, $YXY$, and $Y^{2}X$ are equipollent monomials, and $%
3XY^{2}+3YZ^{2}+3ZX^{2}$ and $%
XY^{2}+YXY+Y^{2}X+YZ^{2}+ZYZ+Z^{2}Y+ZX^{2}+XZX+X^{2}Z$ are equipollent
polynomials.
\end{definition}

\begin{definition}
\label{def3.2}A polynomial is said to be \textit{regular} (or \textit{normal}%
) if it can be expressed as a linear combination of powers of the form 
\begin{equation}
(\alpha X+\beta Y+\gamma Z+\ldots )^{p},\qquad p\in \mathbb{N},\;\alpha
,\beta ,\gamma \in \mathbb{K}.
\end{equation}
\end{definition}

Poincar\'{e} then makes several statements without bothering to prove them.
(They must seem to be obvious to him; note that the same statements are made
in \cite{Po4}, which is an abridged version of \cite{Po1}, where regular
polynomials are called normal.) However, because of the importance of their
implications we shall formulate these statements as a theorem and provide
the reader with a proof which seems to be quite long, but we do not see how
to shorten it (we suspect that because of these claims some authors did not
consider Poincar\'{e}'s proof rigorous).

\begin{theorem}
\label{Theorem3.3}$\mathrm{(i)}$~\textit{A necessary and sufficient
condition for a polynomial to be regular is that if it contains among its
terms a certain monomial then it must contain all monomials equipollent to
that monomial and with the same coefficient}.

$\mathrm{(ii)}$~\textit{Among all polynomials equipollent to a given
polynomial there exists one and only one regular polynomial}$.$
\end{theorem}

Some preparatory work is needed for the proof of this theorem.

Let $\mathcal{P}\equiv \mathcal{P}(x_{1},\ldots ,x_{n})$ denote the \textit{%
commutative} algebra of polynomials in $n$ indeterminates $x_{1},\ldots
,x_{n}$, with coefficients in the field $\mathbb{K}$. For an integer $m\geq
0 $ let $\mathcal{P}^{m}\equiv \mathcal{P}^{m}(x_{1},\ldots ,x_{n})$ denote
the subspace of all homogeneous polynomials of degree $m$. If $(\alpha
)=(\alpha _{1},\ldots ,\alpha _{n})$ is a multi-index of non-negative
integers set $x^{(\alpha )}=x_{1}^{\alpha _{1}}\cdots x_{n}^{\alpha _{n}}$
and $\left| (\alpha )\right| =\sum_{i=1}^{n}\alpha _{i}.$ Then it is clear
that the set $\{x^{(\alpha )}\}$ where $(\alpha )$ ranges over all
multi-indices such that $\left| (\alpha )\right| =m$ forms a basis for $%
\mathcal{P}^{m}$. Let $\mathcal{A}$ $\equiv \mathcal{A}(X_{1,}\ldots ,X_{n})$
denote the \textit{algebra of non-commutative polynomials} in $X_{1},\ldots
,X_{n}$, then $\mathcal{A}$ is obviously graded. Let $\mathcal{A}^{m}\equiv 
\mathcal{A}^{m}(X_{1},\ldots ,X_{n})$ denote the subspace of all homogeneous
elements of $\mathcal{A}$ of degree $m$. Define the \textit{symmetrization
map} (Poincar\'{e} does not formulate this map explicitly, but it is obvious
from the context that he must have it in mind) $\Phi _{m}:\mathcal{P}%
^{m}\rightarrow \mathcal{A}^{m}$ as follows:

For $1\leq j\leq \alpha _{1}$ let $X_{j}^{\prime }=X_{1}$, for $\alpha
_{1}+1\leq j\leq \alpha _{1}+\alpha _{2}$ set $X_{j}^{\prime }=X_{2},$ for $%
\alpha _{1}+\alpha _{2}+1\leq j\leq \alpha _{1}+\alpha _{2}+\alpha _{3}$,
set $X_{j}^{\prime }=X_{3},\ldots $, for $\alpha _{1}+\cdots +\alpha
_{n-1}+1\leq j\leq m$ set $X_{j}^{\prime }=X_{n}$, with the convention that
whenever $\alpha _{i}=0$ then the term $X_{i}$ does not appear. Let 
\begin{equation}
\Phi _{m}(x^{(\alpha )})=\limfunc{Sym}(X^{(\alpha )})=\frac{1}{m!}%
\sum_{\sigma \in \Sigma _{m}}X_{\sigma (1)}^{^{\prime }}\cdots X_{\sigma
(m)}^{\prime }  \label{eq3.2}
\end{equation}
and extend by linearity on all elements of $\mathcal{P}^{m}$ (note that $%
X^{(\alpha )}=X_{1}^{\alpha _{1}}\cdots X_{n}^{\alpha _{n}})$. For example,
with $(\alpha )=(1,2,0,\ldots ,0)$ and $m=3$, then $x^{(\alpha
)}=x_{1}x_{2}^{2}$, $X^{(\alpha )}=X_{1}X_{2}^{2}$, $X_{1}^{\prime }=X_{1}$, 
$X_{2}^{\prime }=X_{2}$, $X_{3}^{\prime }=X_{2}$, and 
\begin{align}
\limfunc{Sym}(X^{(\alpha )})& =\frac{1}{6}(X_{1}^{\prime }X_{2}^{\prime
}X_{3}^{\prime }+X_{1}^{\prime }X_{3}^{\prime }X_{2}^{\prime }+X_{3}^{\prime
}X_{2}^{\prime }X_{1}^{\prime }  \label{eq3.3} \\
& \qquad +X_{2}^{\prime }X_{1}^{\prime }X_{3}^{\prime }+X_{2}^{\prime
}X_{3}^{\prime }X_{1}^{\prime }+X_{3}^{\prime }X_{1}^{\prime }X_{2}^{\prime
})  \notag \\
& =\frac{1}{6}%
(X_{1}X_{2}^{2}+X_{1}X_{2}^{2}+X_{2}^{2}X_{1}+X_{2}X_{1}X_{2}+X_{2}^{2}X_{1}+X_{2}X_{1}X_{2})
\notag \\
& =\frac{1}{3}(X_{1}X_{2}^{2}+X_{2}X_{1}X_{2}+X_{2}^{2}X_{1}).  \notag
\end{align}
Now the dimension of $\mathcal{P}^{m}$ is precisely the number of ways a $n$%
-tuple of integers $(\alpha _{1},\ldots ,\alpha _{n})$ can be chosen so that 
$\left| (\alpha )\right| =m$. If $d_{m}$ denotes this number then a
combinatorial formula gives $d_{m}=\binom{n+m-1}{m}$. Order this set $%
\Lambda _{m}$ of multi-indices following the reverse lexicographic ordering
as follows: 
\begin{eqnarray}
&&\text{``}(\alpha )\prec (\beta )\text{ if for some }k,1\leq k\leq n,\alpha
_{k}<\beta _{k}\text{ }  \label{eq3.4} \\
&&\text{and }\alpha _{k+1}=\beta _{k+1},\ldots ,\alpha _{n}=\beta _{n}, 
\notag \\
&&\text{and }(\alpha )\preceq (\beta )\text{ if either }(\alpha )\prec
(\beta )\text{ or }(\alpha )=(\beta )\text{''.}  \notag
\end{eqnarray}
Then $\preceq $ is obviously a total ordering; for example, $(m,0,\ldots
,0)\prec (m-1,1,0,\ldots ,0)\prec (m-3,3,0,\ldots ,0)$, and $(m,0,\ldots ,0)$
is the first (the least) element and $(0,0,\ldots ,0,m)$ is the last (the
largest) element under this ordering. If $p$ is an element of $\mathcal{P}%
^{m}$ of the form 
\begin{equation*}
p=(c_{1}x_{1}+\cdots +c_{n}x_{n})^{m}\text{, \qquad where }c_{i}\in \mathbb{K%
},1\leq i\leq n,
\end{equation*}
then clearly 
\begin{equation*}
p=\sum_{(\alpha )\in \Lambda _{m}}\binom{m}{\alpha }c^{(\alpha )}x^{(\alpha
)},
\end{equation*}
where the multinomial coefficient $\binom{m}{\alpha }$ is equal to $\frac{m!%
}{\alpha _{1}!\cdots \alpha _{n}!}$, and $c^{(\alpha )}=c_{1}^{\alpha
_{1}}\cdots c_{n}^{\alpha _{n}}$. If we denote by $\tilde{X}^{(\alpha )}$
the image of $x^{(\alpha )}$ by $\Phi _{m}$ (i.e., $\tilde{X}^{(\alpha )}=%
\limfunc{Sym}(X^{\alpha })$), then by linearity 
\begin{equation*}
\Phi _{m}(p)=\sum_{(\alpha )\in \Lambda _{m}}\binom{m}{\alpha }c^{(\alpha )}%
\tilde{X}^{(\alpha )}.
\end{equation*}
On the other hand an easy computation shows that 
\begin{equation}
\binom{m}{\alpha }\tilde{X}^{(\alpha )}=\text{sum of all distinct elements
of }\mathcal{A}^{m}\text{ equipollent to }X^{(\alpha )}\text{.}
\label{eq3.5}
\end{equation}
By expanding $(c_{1}X_{1}+\cdots +c_{n}X_{n})^{m}$ and taking into account
the non-commutativity of the products of the $X_{i}$ we see that 
\begin{eqnarray}
&&(c_{1}X_{1}+\cdots +c_{n}X_{n})^{m}  \label{eq3.6} \\
&&\qquad =\sum_{(\alpha )\in \Lambda _{m}}c^{(\alpha )}\text{(sum of all
distinct elements of }\mathcal{A}^{m}\text{ equipollent to }X^{(\alpha )}%
\text{)}  \notag \\
&&\qquad =\sum_{(\alpha )\in \Lambda _{m}}\binom{m}{\alpha }c^{(\alpha )}%
\tilde{X}^{(\alpha )}\text{.}  \notag \\
&&\text{Thus }\Phi _{m}((c_{1}x_{1}+\cdots
+c_{n}x_{n})^{m})=(c_{1}X_{1}+\cdots +c_{n}X_{n})^{m}.  \notag
\end{eqnarray}
It follows by linearity that the image of a regular polynomial in $\mathcal{P%
}^{m}$ is the regular of $\mathcal{A}^{m}$, obtained by substituting the
variable $x_{i}$ by the variable $X_{i}$; moreover, \textit{regular elements
of }$\mathcal{A}^{m}$\textit{\ are already symmetrized}. Let $\Phi :\mathcal{%
P}\rightarrow \mathcal{A}$ denote the linear map obtained by setting 
\begin{equation}
\Phi (p)=\sum_{m\geq 0}\Phi _{m}(p_{m}),  \label{eq3.7}
\end{equation}
where $p$ is decomposed into homogeneous elements as $p=\sum_{m\geq 0}p_{m}.$

If $\sum_{(\alpha )\in \Lambda _{m}}\lambda _{(\alpha )}x^{(\alpha )}$ is an
arbitrary element of $\mathcal{P}^{m}(\lambda _{(\alpha )}\in \mathbb{K}%
,\forall (\alpha )\in \Lambda _{m})$, then from Eqs.\ (\ref{eq3.2}) and (\ref
{eq3.5}) it follows that 
\begin{multline*}
\Phi _{m}\left( \sum_{(\alpha )\in \Lambda _{m}}\lambda _{(\alpha
)}x^{(\alpha )}\right) =\underset{(\alpha )\in \Lambda _{m}}{\sum }\lambda
_{(\alpha )}\tilde{X}^{(\alpha )} \\
=\sum_{(\alpha )\in \Lambda _{m}}\lambda _{(\alpha )}\frac{1}{\binom{m}{%
\alpha }}\text{ (sum of all distinct elements of }\mathcal{A}^{m}\text{
equipollent to }X^{(\alpha )}\text{).}
\end{multline*}
Since the non-ordered monomials of degree $m$ (i.e., the set of all distinct
elements of $\mathcal{A}^{m}$ equipollent to $X^{(\alpha )}$ for all $%
(\alpha )\in \Lambda _{m}$) form a basis for $\mathcal{A}^{m}$, it follows
that if 
\begin{equation*}
\Phi _{m}\left( \sum_{(\alpha )\in \Lambda _{m}}\lambda _{(\alpha
)}x^{(\alpha )}\right) =0,
\end{equation*}
then $\lambda _{(\alpha )}=0$ for all $(\alpha )\in \Lambda _{m}$. Thus $%
\Phi _{m}$ is one-to-one, and it follows that $\Phi $ is a \textit{%
monomorphism} of vector spaces (but \textit{not} an algebra homomorphism).
Let $\mathcal{R}\equiv \mathcal{R}(X_{1},\ldots ,X_{n})$ denote the subspace
of all regular polynomials of $\mathcal{A}$; then among the consequences of
Theorem \ref{Theorem3.3} one can infer that $\Phi $\textit{\ is a vector
space isomorphism of the vector space }$\mathcal{P}$\textit{\ onto the
vector space }$\mathcal{R}$, \emph{and} \emph{moreover, regular elements in }%
$\mathcal{A}$\emph{\ are already symmetrized}.

\begin{lemma}
\label{lemma3.4}For any positive integer $r$\ there exists an $n$-tuple $%
(c)=(c_{1},\ldots ,c_{n})$\ of positive integers such that 
\begin{equation*}
(c)^{ (\beta )}\geq rc^{ (\alpha )}\text{ whenever }(\alpha )\prec (\beta
)\;\forall (\alpha ),(\beta )\in \Lambda _{m}.
\end{equation*}
\end{lemma}

\begin{proof}
Given $r\in \mathbb{N}^{\ast }$ choose $c_{1}=1$ and define $c_{k}$
inductively by $c_{k}=r(c_{k-1})^{m}$ for $2\leq k\leq n$. Then clearly $%
c_{1}\leq c_{2}\leq \cdots \leq c_{n}$ and $c_{k}=r(c_{k-1})^{m}\geq
rc_{1}^{\alpha _{1}}c_{2}^{\alpha _{2}}\cdots c_{k-1}^{\alpha _{k-1}}$ since 
$\alpha _{1}+\cdots +\alpha _{n}=m$. Hence if $(\alpha )\prec (\beta )$,
i.e., $\alpha _{k}<\beta _{k}$ and $\alpha _{k+1}=\beta _{k+1},\ldots
,\alpha _{n}=\beta _{n}$ then 
\begin{multline*}
c^{(\beta )}\geq c_{k}^{\beta _{k}}c_{k+1}^{\beta _{k+1}}\cdots c_{n}^{\beta
_{n}}\geq c_{k}c_{k}^{\alpha _{k}}c_{k+1}^{\alpha _{k+1}}\cdots
c_{n}^{\alpha _{n}} \\
\geq \left( rc_{1}^{\alpha _{1}}c_{2}^{\alpha _{2}}\cdots c_{k-1}^{\alpha
_{k-1}}\right) \left( c_{k}^{\alpha _{k}}c_{k+1}^{\alpha _{k+1}}\cdots
c_{n}^{\alpha _{n}}\right) =rc^{(\alpha )}.
\end{multline*}
\end{proof}

\begin{lemma}
\label{lemma3.5}\textit{The polynomials} $\tilde{X}^{(\alpha )}\equiv \Phi
_{m}(x^{(\alpha )})\equiv \limfunc{Sym}(X^{\alpha })$\textit{\ are regular
for all }$(\alpha )\in \Lambda _{m}$.
\end{lemma}

\begin{proof}
We prove by induction that for each $(\beta )\in \Lambda _{m}$ the following
statement holds: 
\begin{eqnarray}
&&\text{\emph{``For every} }(\alpha )\preceq (\beta ),\tilde{X}^{(\alpha )}%
\text{\emph{can be expressed as}}  \label{eq3.8} \\
&&\quad \quad \quad \tilde{X}^{(\alpha )}=f_{(\alpha )}^{(\beta )}+\underset{%
(\gamma )\succ (\beta )}{\sum }\lambda _{(\alpha )(\gamma )}^{(\beta )}%
\tilde{X}^{(\gamma )},  \notag \\
&&\text{\emph{where} }f_{(\alpha )}^{(\beta )}\text{\emph{is a regular
element of} }\mathcal{A}^{m}  \notag \\
&&\text{\emph{and the constants} }\lambda _{(\alpha )(\gamma )}^{(\beta )}%
\text{\emph{are rational numbers''.}}  \notag
\end{eqnarray}
First observe that if 
\begin{equation*}
(\beta )=(0,\ldots ,0,\underset{i\text{th slot}}{\underset{\uparrow }{m,}}%
0,\ldots ,0),
\end{equation*}
then $x^{(\beta )}=x_{i}^{m}$ and $\Phi _{m}(x^{(\beta )})=\tilde{X}^{(\beta
)}=X_{i}^{m}$, which is by definition regular. Thus the first element in
this reverse lexicographical ordering is $(\beta )=(m,0,\ldots ,0)$ and the
statement (\ref{eq3.8}) holds trivially with $\tilde{X}^{(m,0,\ldots
,0)}=X_{1}^{m}$, where 
\begin{equation*}
f_{(m,0,\ldots ,0)}^{(m,0,\ldots ,0)}=X_{1}^{m}
\end{equation*}
and 
\begin{equation*}
\lambda _{(m,0,\ldots ,0),(\gamma )}^{(m,0,\ldots ,0)}=0\;\forall (\gamma
)\succ (m,0,\ldots ,0).
\end{equation*}
Now assume the statement holds for $(\beta )$. Let $(\beta ^{\prime })$
denote the immediate successor to $(\beta )$, and consider the element $%
g_{(\beta ^{\prime })}\in \mathcal{A}^{m}$ of the form 
\begin{equation*}
g_{(\beta ^{\prime })}=(c_{1}X_{1}+\cdots +c_{n}X_{n})^{m}=\underset{(\alpha
)\in \Lambda _{m}}{\sum }\binom{m}{(\alpha )}c^{(\alpha )}\tilde{X}^{(\alpha
)},
\end{equation*}
where the $n$-tuple $(c)=(c_{1},\ldots ,c_{n})$ is yet to be determined.
This implies that 
\begin{equation*}
\binom{m}{(\beta ^{\prime })}c^{(\beta ^{\prime })}\tilde{X}^{(\beta
^{\prime })}=g_{(\beta ^{\prime })}-\sum_{(\alpha )\preceq (\beta )}\left( 
\begin{array}{c}
m \\ 
(\alpha )
\end{array}
\right) c^{(\alpha )}\tilde{X}^{(\alpha )}-\sum_{(\gamma )\succ (\beta
^{\prime })}\left( 
\begin{array}{c}
m \\ 
(\gamma )
\end{array}
\right) c^{(\gamma )}\tilde{X}^{(\gamma )}.
\end{equation*}
>From (\ref{eq3.8}) it follows that 
\begin{align}
\left( 
\begin{array}{c}
m \\ 
(\beta ^{\prime })
\end{array}
\right) c^{(\beta ^{\prime })}\tilde{X}^{(\beta ^{\prime })}& =g_{(\beta
^{\prime })}-\sum_{(\alpha )\preceq (\beta )}\left( 
\begin{array}{c}
m \\ 
(\alpha )
\end{array}
\right) c^{(\alpha )}\left( f_{(\alpha )}^{(\beta )}+\sum_{(\gamma )\succ
(\beta )}\lambda _{(\alpha )(\gamma )}^{(\beta )}\tilde{X}^{(\gamma )}\right)
\label{eq3.9} \\
& \qquad -\sum_{(\gamma )\succ (\beta ^{\prime })}\left( 
\begin{array}{c}
m \\ 
(\gamma )
\end{array}
\right) c^{(\gamma )}\tilde{X}^{(\gamma )}  \notag \\
& =\left( g_{(\beta ^{\prime })}-\sum_{(\alpha )\preceq (\beta )}\left( 
\begin{array}{c}
m \\ 
(\alpha )
\end{array}
\right) c^{(\alpha )}f_{(\alpha )}^{(\beta )}\right)  \notag \\
& \qquad -\sum_{(\alpha )\preceq (\beta )}\sum_{(\gamma )\succ (\beta
)}\left( 
\begin{array}{c}
m \\ 
(\alpha )
\end{array}
\right) c^{(\alpha )}\lambda _{(\alpha )(\gamma )}^{(\beta )}\tilde{X}%
^{(\gamma )}  \notag \\
& \qquad -\sum_{(\gamma )\succ (\beta ^{\prime })}\left( 
\begin{array}{c}
m \\ 
(\gamma )
\end{array}
\right) c^{(\gamma )}\tilde{X}^{(\gamma )}\text{.}  \notag
\end{align}
Since $(\beta ^{\prime })$ is right after $(\beta )$, the multi-indices $%
(\gamma )\succ (\beta )$ consist of $(\gamma )=(\beta ^{\prime })$ and $%
(\gamma )\succ (\beta ^{\prime })$. Hence Eq.\ (\ref{eq3.9}) can be written
as 
\begin{multline*}
\left( 
\begin{array}{c}
m \\ 
(\beta ^{\prime })
\end{array}
\right) c^{(\beta ^{\prime })}\tilde{X}^{(\beta ^{\prime })}=\left(
g_{(\beta ^{\prime })}-\sum_{(\alpha )\preceq (\beta )}\left( 
\begin{array}{c}
m \\ 
(\alpha )
\end{array}
\right) c^{(\alpha )}f_{(\alpha )}^{(\beta )}\right) \\
-\left( \sum_{(\alpha )\preceq (\beta )}\left( 
\begin{array}{c}
m \\ 
(\alpha )
\end{array}
\right) c^{(\alpha )}\lambda _{(\alpha )(\beta ^{\prime })}^{(\beta
)}\right) \tilde{X}^{(\beta ^{\prime })} \\
-\sum_{(\gamma )\succ (\beta ^{\prime })}\left\{ \left( \sum_{(\alpha
)\preceq (\beta )}\left( 
\begin{array}{c}
m \\ 
(\alpha )
\end{array}
\right) c^{(\alpha )}\lambda _{(\alpha )(\gamma )}^{(\beta )}\right) +\left( 
\begin{array}{c}
m \\ 
(\gamma )
\end{array}
\right) c^{(\gamma )}\right\} \tilde{X}^{(\gamma )}.
\end{multline*}
This implies that 
\begin{eqnarray}
&&\left[ \left( 
\begin{array}{c}
m \\ 
(\beta ^{\prime })
\end{array}
\right) c^{(\beta ^{\prime })}+\sum_{(\alpha )\preceq (\beta )}\left( 
\begin{array}{c}
m \\ 
(\alpha )
\end{array}
\right) c^{(\alpha )}\lambda _{(\alpha )(\beta ^{\prime })}^{(\beta )}\right]
\tilde{X}^{(\beta ^{\prime })}  \label{eq3.10} \\
&=&\left( g_{(\beta ^{\prime })}-\sum_{(\alpha )\preceq (\beta )}\left( 
\begin{array}{c}
m \\ 
(\alpha )
\end{array}
\right) c^{(\alpha )}f_{(\alpha )}^{(\beta )}\right)  \notag \\
&&-\sum_{(\gamma )\succ (\beta ^{\prime })}\left\{ \left( \sum_{(\alpha
)\preceq (\beta )}\left( 
\begin{array}{c}
m \\ 
(\alpha )
\end{array}
\right) c^{(\alpha )}\lambda _{(\alpha )(\gamma )}^{(\beta )}\right) +\left( 
\begin{array}{c}
m \\ 
(\gamma )
\end{array}
\right) c^{(\gamma )}\right\} \tilde{X}^{(\gamma )}.  \notag
\end{eqnarray}
In Eq.\ (\ref{eq3.10}) we can solve for $\tilde{X}^{(\beta ^{\prime })}$
provided that 
\begin{equation*}
\left( 
\begin{array}{c}
m \\ 
(\beta ^{\prime })
\end{array}
\right) c^{(\beta ^{\prime })}+\sum_{(\alpha )\preceq (\beta )}\left( 
\begin{array}{c}
m \\ 
(\alpha )
\end{array}
\right) c^{(\alpha )}\lambda _{(\alpha )(\beta ^{\prime })}^{(\beta )}
\end{equation*}
is not zero. To insure this we now determine $(c)=(c_{1},\ldots ,c_{n})$ as
in Lemma \ref{lemma3.4} by choosing the integer $r$ such that 
\begin{equation}
r>d_{m}\max_{(\alpha )\preceq (\beta )}\left\{ \left( 
\begin{array}{c}
m \\ 
(\alpha )
\end{array}
\right) \left| \lambda _{(\alpha )(\beta ^{\prime })}^{(\beta )}\right|
\right\} ,  \label{eq3.11}
\end{equation}
where $d_{m}=\binom{m+n-1}{m}$ is the cardinality of $\Lambda _{m}$. Then 
\begin{equation*}
\left( 
\begin{array}{c}
m \\ 
(\beta ^{\prime })
\end{array}
\right) c^{(\beta ^{\prime })}\geq c^{(\beta ^{\prime })}>\sum_{(\alpha
)\preceq (\beta )}\frac{c^{(\beta ^{\prime })}}{d_{m}}\geq \sum_{(\alpha
)\preceq (\beta )}\frac{rc^{(\alpha )}}{d_{m}}>\sum_{(\alpha )\preceq (\beta
)}\left( 
\begin{array}{c}
m \\ 
(\alpha )
\end{array}
\right) \left| \lambda _{(\alpha )(\beta ^{\prime })}^{(\beta )}\right|
c^{(\alpha )}.
\end{equation*}
It follows that 
\begin{multline*}
\left( 
\begin{array}{c}
m \\ 
(\beta ^{\prime })
\end{array}
\right) c^{(\beta ^{\prime })}+\sum_{(\alpha )\preceq (\beta )}\left( 
\begin{array}{c}
m \\ 
(\alpha )
\end{array}
\right) c^{(\alpha )}\lambda _{(\alpha )(\beta ^{\prime })}^{(\beta )} \\
\geq \left( 
\begin{array}{c}
m \\ 
(\beta ^{\prime })
\end{array}
\right) c^{(\beta ^{\prime })}-\sum_{(\alpha )\preceq (\beta )}\left( 
\begin{array}{c}
m \\ 
(\alpha )
\end{array}
\right) c^{(\alpha )}\left| \lambda _{(\alpha )(\beta ^{\prime })}^{(\beta
)}\right| >0.
\end{multline*}
Thus we have shown that the $n$-tuple $(c)=(c_{1},\ldots ,c_{n})$ of
positive integers can be chosen so that the coefficient of $\tilde{X}%
^{(\beta ^{\prime })}$ in Eq.\ (\ref{eq3.10}) is a positive rational number,
and the coefficients of $\tilde{X}^{(\gamma )}$ in the sum $\sum_{(\gamma
)\succ (\beta ^{\prime })}$ are rational numbers. Obviously, 
\begin{equation*}
g_{(\beta ^{\prime })}-\sum_{(\alpha )\preceq (\beta )}\binom{m}{(\alpha )}%
c^{(\alpha )}f_{(\alpha )}^{(\beta )}
\end{equation*}
is a regular polynomial, so by dividing both sides of Eq.\ (\ref{eq3.10}) by
the coefficient of $\tilde{X}^{(\beta ^{\prime })}$ we can write 
\begin{equation}
\tilde{X}^{(\beta ^{\prime })}=f_{(\beta ^{\prime })}^{(\beta ^{\prime
})}+\sum_{(\gamma )\succ (\beta ^{\prime })}\lambda _{(\beta ^{\prime
})(\gamma )}^{(\beta ^{\prime })}\tilde{X}^{(\gamma )},  \label{eq3.12}
\end{equation}
where $f_{(\beta ^{\prime })}^{(\beta ^{\prime })}$ is regular and the
constants $\lambda _{(\beta ^{\prime })(\gamma )}^{(\beta ^{\prime })}$ are
rational. For $(\alpha )\prec (\beta ^{\prime })$, i.e., $\alpha \preceq
(\beta )$ Eq.\ (\ref{eq3.8}) can be written as 
\begin{align}
\tilde{X}^{(\alpha )}& =f_{(\alpha )}^{(\beta )}+\lambda _{(\alpha )(\beta
^{\prime })}^{(\beta )}\tilde{X}^{(\beta ^{\prime })}+\sum_{(\gamma )\succ
(\beta ^{\prime })}\lambda _{(\alpha )(\gamma )}^{(\beta )}\tilde{X}%
^{(\gamma )}  \label{eq3.13} \\
& =\left( f_{(\alpha )}^{(\beta )}+\lambda _{(\alpha )(\beta ^{\prime
})}^{(\beta )}f_{(\beta ^{\prime })}^{(\beta ^{\prime })}\right)
+\sum_{(\gamma )\succ (\beta ^{\prime })}\left( \lambda _{(\alpha )(\beta
^{\prime })}^{(\beta )}\lambda _{(\beta ^{\prime })(\gamma )}^{(\beta
^{\prime })}+\lambda _{(\alpha )(\gamma )}^{(\beta )}\right) \tilde{X}%
^{(\gamma )}.  \notag
\end{align}
Set $f_{(\alpha )}^{(\beta ^{\prime })}=f_{(\alpha )}^{(\beta )}+\lambda
_{(\alpha )(\beta ^{\prime })}^{(\beta )}f_{(\beta ^{\prime })}^{(\beta
^{\prime })}$ and $\lambda _{(\alpha )(\gamma )}^{(\beta ^{\prime
})}=\lambda _{(\alpha )(\beta ^{\prime })}^{(\beta ^{\prime })}\lambda
_{(\beta ^{\prime })(\gamma )}^{(\beta ^{\prime })}+\lambda _{(\alpha
)(\gamma )}^{(\beta )}$, then it follows from Eqs.\ (\ref{eq3.12}) and (\ref
{eq3.13}) that for all $(\alpha )\preceq (\beta ^{\prime })$, 
\begin{equation*}
\tilde{X}^{(\alpha )}=f_{(\alpha )}^{(\beta ^{\prime })}+\sum_{(\gamma
)\succ (\beta ^{\prime })}\lambda _{(\alpha )(\gamma )}^{(\beta ^{\prime })}%
\tilde{X}^{(\gamma )},
\end{equation*}
where $f_{(\alpha )}^{(\beta ^{\prime })}$ is obviously regular and the
coefficients $\lambda _{(\alpha )(\gamma )}^{(\beta ^{\prime })}$ are
obviously rational. Hence we have completed the induction. Now for the proof
of the lemma in the statement (\ref{eq3.8}), choose $(\beta )=(\beta )_{\max
}=(0,\ldots ,0,m)$ to be the last element of $\Lambda _{m}$. Then (\ref
{eq3.8}) reads:

``For every $(\alpha )\preceq (\beta )_{\max }$ $\tilde{X}^{(\alpha
)}=f_{(\alpha )}^{(\beta )_{\max }}$, where $f_{(\alpha )}^{(\beta )_{\max
}} $ is a regular element of $\mathcal{A}^{m}$''. This is exactly what the
lemma affirms.
\end{proof}

>From the fact that $\Phi _{m}:\mathcal{P}^{m}\rightarrow \mathcal{A}^{m}$ is
a monomorphism it follows that the system $\{\tilde{X}^{(\alpha )},(\alpha
)\in \Lambda _{m}\}$ is linearly independent. Therefore if $\mathcal{R}^{m}$
denotes the subspace of $\mathcal{R}$ of all homogeneous non-commutative
regular polynomials of degree $m$ in the indeterminates $X_{1},\ldots ,X_{n}$%
, then Lemma \ref{lemma3.5} and Eq.\ (\ref{eq3.6}) imply that \textit{the
system }$\{\tilde{X}^{\alpha },(\alpha )\in \Lambda _{m}\}$\textit{\ forms a
basis for }$\mathcal{R}^{m}$. It follows immediately from the discussion
preceding Lemma \ref{lemma3.4} that $\Phi _{m}:\mathcal{P}^{m}\rightarrow 
\mathcal{R}^{m}$ is an isomorphism and hence, $\Phi $\textit{\ is an
isomorphism of }$\mathcal{P}$\textit{\ onto }$\mathcal{R}$ (clearly from
Eq.\ (\ref{eq3.6}) $\mathcal{R}^{m}\subset \Phi _{m}(\mathcal{P}^{m})$,
Lemma \ref{lemma3.5} shows that $\mathcal{R}^{m}=\Phi _{m}(\mathcal{P}^{m})$%
). Now each $\tilde{X}^{(\alpha )}$, being a regular homogeneous polynomial
of degree $m$, is therefore a linear combination of polynomials of the form $%
(c_{1}X_{1}+\cdots +c_{n}X_{n})^{m}$. Let $S$ be the set of such
polynomials, then $S$ is a finite set of vectors spanning the vector space $%
\mathcal{R}^{m}$. From a general fact in linear algebra (see, for example, 
\cite[Corollary 2, p.~44]{H-K}) we can deduce the following.

\begin{corollary}
\label{corollary3.6}$\mathrm{(i)}$ \textit{The vector space }$\mathcal{R}%
^{m} $\textit{\ admits a basis consisting of vectors of the form} 
\begin{equation}
f_{i}=(c_{1}^{i}X_{1}+c_{2}^{i}X_{2}+\cdots +c_{n}^{i}X_{n})^{m},\qquad
1\leq i\leq d_{m},c_{j}^{i}\in \mathbb{K}.  \label{eq3.14}
\end{equation}
$\mathrm{(ii)}$~\textit{The same conclusion holds with }$\mathcal{P}^{m}$%
\textit{\ replacing }$\mathcal{R}^{m}$\textit{\ and }$x_{j},1\leq j\leq n$%
\textit{\ replacing }$X_{j\text{.}}$
\end{corollary}

\textit{Proof of Theorem }\textrm{\ref{Theorem3.3}}. Clearly since $\mathcal{%
A}$ is a graded algebra, it suffices to prove the theorem for $\mathcal{A}%
^{m},m\geq 0$.

(i)~If $f$ is a regular element of $\mathcal{A}^{m}$ then since $\{\tilde{X}%
^{(\alpha )},(\alpha )\in \Lambda _{m}\}$ forms a basis for $\mathcal{R}^{m}$%
, 
\begin{equation*}
f=\sum_{(\alpha )\in \Lambda _{m}}\lambda _{(\alpha )}\tilde{X}^{(\alpha
)},\qquad \lambda _{(\alpha )}\in \mathbb{K}.
\end{equation*}
Thus if $\lambda _{(\alpha )}\neq 0$ then since $\tilde{X}^{(\alpha )}$
contains all distinct monomials equipollent to $X^{(\alpha )}$, therefore $f$
contains all monomials equipollent to $X^{(\alpha )}$ with the same
coefficient 
\begin{equation*}
\frac{1}{\binom{m}{(\alpha )}}\lambda _{(\alpha )}.
\end{equation*}
Conversely, if $f$ is a polynomial in $\mathcal{A}^{m}$ which contains all
monomials equipollent to a fixed monomial, which we may assume without loss
of generality to be $X^{(\alpha )}$, with the same coefficient, then $f$
must contain $\mu _{(\alpha )}\tilde{X}^{(\alpha )}$ with $\mu _{(\alpha )}$
a non-zero constant. Hence $f$ is of the form 
\begin{equation*}
\sum_{(\alpha )\in \Lambda _{m}}\lambda _{(\alpha )}\tilde{X}^{(\alpha
)},\qquad \lambda _{(\alpha )}\in \mathbb{K},
\end{equation*}
and therefore is a regular polynomial.

(ii)~First observe that if $X_{i_{1}}\cdots X_{i_{m}},1\leq i_{j}\leq
n,1\leq j\leq m$, is a monomial in $\mathcal{A}^{m}$, then it is equipollent
to a unique monomial $X^{(\alpha )}=X_{1}^{\alpha _{1}}\cdots X_{n}^{\alpha
_{n}}$ for some $(\alpha )\in \Lambda _{m}$. Then from the definition of $%
\tilde{X}^{(\alpha )}\equiv \limfunc{Sym}(X^{\alpha })$ and part (i)~of this
theorem, $\tilde{X}^{(\alpha )}$ is the unique regular polynomial of $%
\mathcal{A}^{m}$ that is equipollent to $X^{(\alpha )}$, and hence to $%
X_{i_{1}\ldots i_{m}}$ (for example, $X_{2}X_{1}X_{2}$ is equipollent to $%
X^{(1,2,0,\ldots ,0)}=X_{1}X_{2}^{2}$, which from Eq.\ (\ref{eq3.3}) is
equipollent to the regular polynomial $\tilde{X}^{(1,2,0,\ldots ,0)}=\frac{1%
}{3}(X_{1}X_{2}^{2}+X_{2}X_{1}X_{2}+X_{2}^{2}X_{1})=\frac{1}{6}[%
(X_{1}+X_{2})^{3}+(X_{1}-X_{2})^{3}-2X_{1}^{3}])$. Let 
\begin{equation*}
p=\sum_{i_{1},\ldots ,i_{m}}\lambda _{_{i_{1\cdots }i_{m}}}X_{i_{1}}\cdots
X_{i_{m}},
\end{equation*}
where the sum is over all distinct non-commutative homogeneous monomials of
degree $m$ and the coefficients $\lambda _{i_{1}\ldots i_{m}}$ are uniquely
determined. Then since each $X_{i_{1}}\cdots X_{i_{m}}$ is equipollent to a
unique regular polynomial $\tilde{X}^{(\alpha )}$ for some $(\alpha )\in
\Lambda _{m}$, $p$ is equipollent to the unique regular polynomial 
\begin{equation*}
\sum_{(\alpha )\in \Lambda _{m}}\mu _{(\alpha )}\tilde{X}^{(\alpha )},
\end{equation*}
where $\mu _{(\alpha )}$ is the sum of all $\lambda _{i_{1}\ldots i_{m}}$
for which $X_{i_{1}\ldots i_{m}}$ is equipollent to $\tilde{X}^{(\alpha )}$.%
%TCIMACRO{
%\TeXButton{End Proof}{\endproof%
%}}%
%BeginExpansion
\endproof%
%
%EndExpansion

\begin{remark}
\textrm{\label{remark3.7}}It follows from Corollary \ref{corollary3.6}(i)
and Eq.\ (\ref{eq3.6}) that a polynomial $p$ of $\mathcal{A}$ is regular if
and only if $\limfunc{Sym}(p)=p$. Now define (as in \cite[5.6.1]{Go}) a
polynomial 
\begin{equation*}
p=\sum_{i_{1}\ldots i_{m}}\lambda _{i_{1}\ldots i_{m}}X_{i_{1}}\cdots
X_{i_{m}},\qquad 1\leq i_{j}\leq n,1\leq j\leq m,
\end{equation*}
to be \textit{symmetric} if all its coefficients $\lambda _{i_{1}\ldots
i_{m}}$ are symmetric, i.e., $\lambda _{i_{\sigma (1)}\ldots i_{\sigma
(m)}}=\lambda _{i_{1}\ldots i_{m}}$ for all $\sigma \in \Sigma _{m}$. Then
Theorem \ref{Theorem3.3}(i) implies that a polynomial in $\mathcal{A}$ is
regular if and only it is symmetric, since a monomial is equipollent to $%
X_{i_{1}}\cdots X_{i_{m}}$, if and only if it is of the form $X_{i_{\sigma
(1)}}\cdots X_{i_{\sigma (m)}}$ for some $\sigma \in \Sigma _{m}$. Thus in
this context regular is synonymous with symmetric, and this is probably what 
\cite{Bo} must have had in mind when he affirmed that Poincar\'{e} gave a
proof of algebraic nature that the associative algebra generated by the $%
X_{i},1\leq i\leq n$, has as basis certain symmetric functions in $X_{i}$.
In fact, in \cite[5.6.1]{Go}, for example, this fact is used to define the
vector space isomorphism $\beta :S(\frak{g})\rightarrow \mathcal{U}(\frak{g}%
) $ of the symmetric algebra of polynomial functions on the dual $\frak{g}%
\ast $ of the Lie algebra $\frak{g}$ onto the universal enveloping algebra $%
\frak{g}$. Obviously, $\mathcal{P}$ is isomorphic to $S(\frak{g})$ and as we
shall see $\mathcal{R}$ is isomorphic to $\mathcal{U}(\frak{g})$; thus the
map $\beta $ is basically $\Phi $.
\end{remark}

Now let us return to \cite{Po1}. Let $X_{1},\ldots ,X_{r}$ be elementary
operators (i.e., infinitesimal transformations) which form a basis for a Lie
algebra $\mathcal{L}$. Define the Lie bracket as 
\begin{equation}
\lbrack X,Y]=XY-YX;\qquad X,Y\in \mathcal{L}\text{.}  \label{eq3.15}
\end{equation}
Two polynomials in $\mathcal{A}$ are said to be \emph{equivalent} if one can
be reduced to the other in taking into account relation (\ref{eq3.15})$.$

For example, the product $P(XY-YX-[X,Y])Q$ as defined in Eq.\ (\ref{eq3.1})
(where the first and the last factors $P$ and $Q$ are two arbitrary
monomials in $\mathcal{A}$) is equivalent to zero, and obviously so are
linear combinations of products of that form (i.e., $P$ and $Q$ may be taken
to be polynomials). Products of the form (\ref{eq3.1}) are called \textit{%
trinomial products}.

The difference of two monomials which differ only by the order of two
consecutive factors is equivalent to a polynomial of lesser degree. Indeed,
let $X$ and $Y$ be those two consecutive factors. Then our monomials are
written as 
\begin{equation*}
PXYQ\text{ and }PYXQ\text{,}
\end{equation*}
$P$ and $Q$ being two arbitrary monomials, and their difference 
\begin{equation*}
P(XY-YX)Q
\end{equation*}
is equivalent to $P[X,Y]Q$, which has degree one less, since $[X,Y]$ is of
first degree, while $XY$ and $YX$ are of second degree.

Now let $M$ and $M^{\prime }$ be two arbitrary equipollent monomials; that
is, they only differ by the order of their factors. One can find a sequence
of monomials 
\begin{equation*}
M,M_{1},M_{2},\ldots ,M_{p},M^{\prime },
\end{equation*}
in which the first and the last terms are the given monomials and any term
in the sequence differs only from the preceding by the order of two
consecutive factors. The difference $M-M^{\prime }$, which is the sum of the
differences $M-M_{1},M_{1}-M_{2},\ldots ,M_{p}-M^{\prime }$, is therefore
again equivalent to a polynomial of lesser degree.

More generally, the difference of two equipollent polynomials is equivalent
to a polynomial of lesser degree. We now claim the following.

\begin{theorem}
\label{theorem3.8}\textit{In the algebra }$\mathcal{A}$\textit{\ any
arbitrary polynomial is equivalent to a unique regular polynomial}.
\end{theorem}

%TCIMACRO{
%\TeXButton{Proof}{\proof%
%}}%
%BeginExpansion
\proof%
%
%EndExpansion
First let us show that this equivalence relation is \textit{additive}, i.e.,
if $p\sim p^{\prime }$ \footnote{%
Poincar\'{e} used the symbol $=$ to denote this equivalence relation. To
avoid confusion we adopt the more conventional symbol $\sim $.} and $q\sim
q^{\prime }$ then $p+q\sim p^{\prime }+q^{\prime }$. This is obvious since
this is equivalent to $p-p^{\prime }\sim 0$ and $q-q^{\prime }\sim 0$, and
hence $(p+q)-(p^{\prime }+q^{\prime })=(p-p^{\prime })+(q-q^{\prime })\sim
0+0=0$.

Now let $P_{n}$ be an arbitrary polynomial of degree $n$; then by Theorem 
\ref{Theorem3.3}(ii) $P_{n}$ is equipollent to a unique regular polynomial $%
P_{n}^{\prime }$ of the same degree $n$, and by the remark preceding this
theorem, $P_{n}-P_{n}^{\prime }$ is equivalent to a polynomial $P_{n-1}$ of
lesser degree (which we may assume, without loss of generality, of degree $%
n-1$). Therefore, $P_{n}\sim P_{n}^{\prime }+P_{n-1}$, and $P_{n-1}$ is in
turn equipollent to a regular polynomial $P_{n-1}^{\prime }$, and hence 
\begin{equation*}
P_{n}\sim P_{n}^{\prime }+P_{n-1}=P_{n}^{\prime }+P_{n-1}^{\prime
}+(P_{n-1}-P_{n-1}^{\prime })\sim P_{n}^{\prime }+P_{n-1}^{\prime
}+(P_{n-2}),\ldots ,%
%TCIMACRO{
%\TeXButton{TeX field}{\notag%
%}}%
%BeginExpansion
\notag%
%
%EndExpansion
\end{equation*}
and so on; one finally arrives to a polynomial of degree zero which is
obviously regular. Thus one can conclude that 
\begin{equation*}
P_{n}\sim P_{n}^{\prime }+P_{n-1}^{\prime }+P_{n-2}^{\prime }+\cdots ,
\end{equation*}
where the second member is a regular polynomial. We therefore have a means
to reduce any polynomial to a regular polynomial by making use of the
relations (\ref{eq3.15}). It remains to find out if this reduction can be
done uniquely.

Since both the equivalence relation $\sim $ and the notion of regular
polynomials are additive, this problem is equivalent to the following:

\textit{Can a non-identically zero regular polynomial be equivalent to zero}%
? Or equivalently, \textit{can we find a sum of trinomial products of the
form }($\ref{eq3.1}$)\textit{\ which is a non-identically zero regular
polynomial}? All sums of such products are indeed equivalent to zero and
vice-versa. If we define a \textit{regular sum }to be a sum of trinomial
products of the form (\ref{eq3.1}) which is also a regular polynomial then
the answer (\textit{negative}) to this question (and hence to the question
above regarding uniqueness) can be stated as follows:

\begin{lemma}
\label{lemma3.9}\textit{Every regular sum is identically zero}.
\end{lemma}

\noindent \textit{Proof of the lemma.} The degree of a trinomial product (%
\ref{eq3.1}) is clearly $d^{0}(P)+d^{0}(Q)+2$. Thus we call the degree of a
sum $S$ of trinomial products the highest of all the degrees of the products
in $S$ \textit{even though as we shall see when }$S$\textit{\ is a regular
sum the terms of highest degree in these different products mutually cancel
each other.}

The trinomial product (\ref{eq3.1}) can be considered as the sum of two
products, the \textit{binomial product} 
\begin{equation}
P(XY-YX)Q,  \label{eq3.16}
\end{equation}
where we call $PXYQ$ the \textit{positive monomial} and $-PYXQ$ the \textit{%
negative monomial}; and the product 
\begin{equation}
-P[X,Y]Q,  \label{eq3.17}
\end{equation}
which we call the \textit{complementary product}.

Thus if $S$ is an arbitrary sum of trinomial products of degree $p$ and of
degree $<p$ then we can write 
\begin{equation}
S=S_{p}-T_{p}+S_{p-1}-T_{p-1}+\cdots +S_{k}-T_{k}+\cdots +S_{2}-T_{2},
\label{eq3.18}
\end{equation}
where $S_{k},$ $2\leq k\leq p$, is a sum of homogeneous binomial products of
degree $k,$ whereas $-T_{k}$ is the sum of the corresponding complementary
products. First observe that if $S$ is a regular sum then every homogeneous
component of $S$ is also regular since regularity is graded; in particular $%
\mathit{S}_{p}$\textit{\ is regular}. Since $S_{p}$ is a sum of binomial
products of degree $p$ of the form $PXYQ-PYXQ$ and since equipollence is an
additive equivalence relation it follows immediately that $S_{p}$ is
equipollent to zero. But zero is a regular polynomial and Theorem \ref
{Theorem3.3}(ii) implies that two regular polynomials cannot be equipollent
without being identical, and therefore $S_{p}$ must be identically zero.

\begin{remark}
\textrm{\label{remark3.10}}From the discussion above it follows that \textit{%
the degree of a regular sum as we defined it is actually at least one more} 
\textit{than the classical degree of a polynomial in }$\mathcal{A}$.
\end{remark}

Thus in particular when $S$ is a regular sum of degree 3 (actual degree 2)
then 
\begin{equation}
S=S_{3}-T_{3}+S_{2}-T_{2}.  \label{eq3.19}
\end{equation}

Since $S_{3}$ is homogeneous of degree $3$, a typical binomial product of $%
S_{3}$ must be of the form 
\begin{equation*}
(XY-YX)Z\text{ \quad or\quad }Z(XY-YX).
\end{equation*}

Since $S_{3}$ is regular, hence symmetric, Theorem \ref{Theorem3.3} (i)
implies that if the binomial product $(XY-YX)Z=XYZ-YXZ$ occurs in $S_{3}$,
all six monomials equipollent to $XYZ$ (resp. $-YXZ$) must occur in $S_{3}$
with the same coefficient.

Thus $S_{3}$ must be a sum of terms of the form 
\begin{equation}
\sum (XY-YX)Z-\sum Z(XY-YX),  \label{eq3.20}
\end{equation}
where the sign $\sum $ means that one must sum over the term which is
explicitly expressed under the sign and the other two terms obtained by
cyclically permuting the three letters $X$, $Y$, $Z$. Note that one can
verify directly from Eq.\ (\ref{eq3.20}) that $S_{3}$ is identically zero.
It follows from Eq.\ (\ref{eq3.20}) that the sum of the complementary
products $-T_{3}$ contains terms of the form 
\begin{equation}
-\left( \sum [X,Y]Z-\sum Z[X,Y]\right)  \label{eq3.21}
\end{equation}
Since $S_{2}-T_{3}$ is homogeneous of degree two and $S$ is regular, it
follows that $S_{2}-T_{3}$ is also regular, and hence symmetric. Since it
contains $-[X,Y]Z+Z[X,Y]$, Theorem \ref{Theorem3.3}(i) again implies that it
must contain permutations of these terms with the same coefficients, i.e., 
\begin{multline*}
-[X,Y]Z-[Y,X]Z-[Z,Y]X-[X,Z]Y-[Z,X]Y-[Y,Z]X \\
+Z[X,Y]+X[Z,Y]+Y[Z,X]+Z[Y,X]+Y[X,Z]+X[Y,Z],
\end{multline*}
which can be regrouped in the following form: 
\begin{equation}
-\left( \sum [X,Y]Z-\sum Z[X,Y]\right) +\left( \sum [X,Y]Z-\sum
Z[X,Y]\right) ,  \label{eq.3.22}
\end{equation}
using the fact that the bracket product $[\;,\;]$ is anti-symmetric. From
Eqs.\ (\ref{eq3.21}) and (\ref{eq.3.22}) it follows that 
\begin{equation}
S_{2}=\sum [X,Y]Z-\sum Z[X,Y]  \label{eq3.23}
\end{equation}
and $S_{2}-T_{3}=0$. Since $S_{2}$ is a sum of terms of the form $WZ-ZW$
with $W=[X,Y]$ it follows that the complementary polynomial is of degree one
and is a sum of terms of the form $[W,Z].$ Thus we have 
\begin{equation*}
T_{2}=\sum [[X,Y],Z],
\end{equation*}
where $\sum $ has the same meaning as above. Thus $T_{2}$ is a sum of terms
of the form 
\begin{equation}
\lbrack \lbrack X,Y],Z]+[[Y,Z],X]+[[Z,X],Y].  \label{eq3.24}
\end{equation}

It follows from Eq.\ (\ref{eq2.1}) that $T_{2}$ is a polynomial of first
degree which is obviously symmetric, and hence regular. Therefore, if $T_{2}$
is not identically zero, the sum $S$ would be a regular polynomial which is
not identically zero.

Therefore, in order that a polynomial can be reduced in a unique fashion to
a regular polynomial, it is necessary that the expression (\ref{eq3.24}) is
identically zero. But one recognizes there the Jacobi identities which play
such an important role in Lie theory. It remains to show that this condition
is sufficient.

At this juncture, it is important to make the following remark:

\textit{It follows from Remark }$\ref{remark3.10}$\textit{\ that we have
actually proved that every polynomial of degree }$0$\textit{, }$1$\textit{,
or }$2$\textit{\ is equivalent to a unique regular polynomial of degree }$0$%
\textit{, }$1$\textit{, or }$2$, \textit{respectively.}

Now by induction suppose that the lemma has been proven for regular sums of
degree $1,2,\ldots ,p-1$ and propose to extend it to regular sums of degree $%
p.$

Thus, let $S=S_{p}-T_{p}+S_{p-1}-T_{p-1}+\cdots $ be a sum of trinomial
products. Let us call $S_{p}-T_{p}$ the \textit{head} (or leading terms) of
the sum $S$. We say that a sum of trinomial products form a \textit{chain}
if the negative monomial of each product is equal and of opposite sign of
the positive monomial of the product that follows. The positive monomial of
the first product and the negative monomial of the last one are called
extreme monomials of the chain. Examples of chains: 
\begin{multline*}
C_{1}\colon XZ(XY)W-XZ(YX)W-XZ[X,Y]W+X(ZY)XW-X(YZ)XW \\
-X[Z,Y]XW+XY(ZX)W-XY(XZ)W-XY[Z,X]W.
\end{multline*}
\begin{multline*}
C_{2}\colon XZ(XY)W-XZ(YX)W-XZ[X,Y]W+XZY(XW)-XZY(WX) \\
-XZY[X,W]+X(ZY)WX-X(YZ)WX-X[Z,Y]WX+XYZ(WX) \\
-XYZ(XW)-XYZ[W,X]+XY(ZX)W-XY(XZ)W-XY[Z,X]W.
\end{multline*}

\begin{remark}
\label{remark3.11}\textit{It results from the definition that all positive
monomials }(\textit{and hence, all negative monomials})\textit{\ of the same
chain can only differ by the order of their factors}.
\end{remark}

A chain is said to be \textit{closed} if its extreme monomials are equal and
of opposite sign. If $S_{p}-T_{p}$ is a closed chain of trinomial products
it is clear that $S_{p}$ is identically zero since the positive and negative
monomials cancel each other two by two.

We have seen that if $S$ is a regular sum, $S_{p}$ is identically zero. It
follows therefore that \textit{the head of a regular sum must always consist
of one or more closed chains}.

If two chains have the same extreme monomials, then their difference is a
closed chain. For example, 
\begin{multline*}
C_{1}-C_{2}\colon XYZ(XW)-XYZ(WX)-XYZ[X,W] \\
+X(YZ)WX-X(ZY)WX-X[Y,Z]WX+XZYWX \\
-XYZXW+XZY[X,W]-X[Z,Y]XW.
\end{multline*}
We shall use this remark to show that a closed chain can always be
decomposed in many ways into two or more closed chains. An arbitrary closed
chain can be in many ways regarded as the difference of two chains $C$ and $%
C^{\prime }$ having the same extreme monomials. Let $C^{\prime \prime }$ be
a third chain having the same extreme monomials, then the chain $C-C^{\prime
}$ is then decomposed into two other closed chains $C-C^{\prime \prime }$
and $C^{\prime \prime }-C^{\prime }$.

Now remark first that if a regular sum of degree $p$ is identically zero, it
must be the same for all regular sums of degree $p$ \textit{which have the
same head}. The difference of these two sums will be indeed a regular sum of
degree $p-1$ which will be identically zero according to our inductive
hypothesis. \textit{Therefore it suffices for us to form all closed chains
of degree }$p$\textit{\ and prove that each one of them can be considered as
the head of an identically zero regular sum}. Indeed, each regular sum $S$
of order $p$ has as head one or more of those closed chains. Let $S^{\prime
} $ be one of those closed chains, then if we show that there exists an
identically zero regular sum having $S^{\prime }$ as head, it follows
immediately from the remark above that $S$ must be identically zero. Thus by
induction we suppose this statement holds for all closed chains of degree $%
\leq p-1$ and we will show that it is true for all closed chains of degree $%
p $.

To establish this assertion, we are going to decompose the closed chain in
question into several closed chains. It is clear that it suffices to prove
the proposition for each component.

A chain is called \textit{simple of the first kind} if the first factor of
all of its monomials either positive or negative is everywhere the same. A
chain is called \textit{simple of the second kind} if the last factor of all
of its monomials either positive or negative is everywhere the same.
Moreover, a simple chain can be either closed or \textit{open} (not closed).

Since $p$ is larger than three, it is clear that every closed chain can be
regarded as the sum of a certain number of simple chains, alternatively of
the first and second kinds or vice-versa.

Thus let $S$ be a closed chain, $C_{1},C_{2},\ldots ,C_{n}$ be simple chains
of the first kind, $C_{1}^{\prime },C_{2}^{\prime },\ldots ,C_{n}^{\prime }$
be simple chains of the second kind, such that 
\begin{equation*}
S=C_{1}+C_{1}^{\prime }+C_{2}+C_{2}^{\prime }+\cdots +C_{n}+C_{n}^{\prime },
\end{equation*}
the extreme negative monomial of each chain being, of course, equal and of
opposite sign to the extreme positive monomial of the next chain, and the
extreme negative monomial of $C_{n}^{\prime }$ being equal and of opposite
sign to the extreme positive monomial of $C_{1}$. Note that, a priori, $%
C_{1} $ or $C_{n}^{\prime }$, can be the zero chain, for example, if $S$
starts with $XYQ-YXQ-\cdots $, where $d^{0}(Q)>1$; then $C_{1}=0$, but then
we can consider $C_{1}$ as the zero simple closed chain of the form $%
XYQ-\cdots -XYQ $, and similarly for $C_{n}^{\prime }$.

Let $X$ be the first factor of all the monomials of $C_{1}$, $Z$ the last
factor of all the monomials of $C_{1}^{\prime }$, $Y$ the first factor of
all the monomials of $C_{2}$, and $T$ the last factor of all the monomials
of $C_{2}^{\prime }$ (we do not exclude the case where two of the operators $%
X$, $Y$, $Z$, $T$ are identical).

Let $C^{\prime \prime }$ be a simple chain of the second kind having its
extreme positive monomial equal and of opposite sign to the extreme negative
monomial of $C_{2}^{\prime }$, and in which all monomials have the last
factor equal to $T$, and moreover, the extreme negative monomial has $X$ as
its first factor.

Let $C^{\prime \prime \prime }$ be a simple chain of the first kind such
that all monomials in it have $X$ as the first factor, and moreover, the
extreme monomials are respectively equal and of opposite signs to the
extreme negative monomial of $C^{\prime \prime }$ and to the extreme
positive monomial of $C_{1}$.

Schematically we have the following diagram: 
\begin{align}
& \underset{C_{1}}{\underbrace{X\square \cdots -X\square Z}}+\underset{%
C_{1}^{\prime }}{\underbrace{X\square Z\cdots -Y\square Z}}+\underset{C_{2}}{%
\underbrace{Y\square Z\cdots -Y\square T}}  \label{eq3.25} \\
& \hspace{0.5in}+\underset{C_{2}^{\prime }}{\underbrace{Y\square T\cdots
-\square T}}+\underset{C^{\prime \prime }}{\underbrace{\square T\cdots
-X\square T}}+\underset{C^{\prime \prime \prime }}{\underbrace{X\square
T\cdots -X\square }},  \notag
\end{align}
where each box $\square $ represents certain unspecified monomial which does
not have any effect on our discussion.

Thus the closed chain $S$ is decomposed into a sum of two closed chains as $%
S^{\prime }+S^{\prime \prime }$, where 
\begin{eqnarray*}
S^{\prime } &=&(C^{\prime \prime \prime }+C_{1})+C_{1}^{\prime
}+C_{2}+(C_{2}^{\prime }+C^{\prime \prime }), \\
S^{\prime \prime } &=&-C^{\prime \prime }+C_{3}+\cdots +C_{n}-C^{\prime
\prime \prime }.
\end{eqnarray*}
The closed chain $S^{\prime }$ contains only four simple chains, since $%
(C^{\prime \prime \prime }+C_{1})$ and $(C_{2}^{\prime }+C^{^{\prime \prime
}})$ are simple chains; $S^{\prime \prime }$ contains two simple chains less
than $S$. Continuing this scheme we end up decomposing $S$ into closed
components which consist of only four simple chains. \textit{Thus it
suffices to consider the case of closed chains }$S$\textit{\ formed by four
simple chains as, for example, the form }$S^{\prime }$.

Therefore, it follows from (\ref{eq3.25}) that the extreme positive
monomials of the four chains that form $S^{\prime }$ have respectively for
first and last factors: 
\begin{equation*}
\begin{tabular}{ll}
for $C^{\prime \prime \prime }+C_{1}$ & $X$ and $T,$ \\ 
for $C_{1}^{\prime }$ & $X$ and $Z,$ \\ 
for $C_{2}$ & $Y$ and $Z,$ \\ 
for $C_{2}^{\prime }+C^{\prime \prime }$ & $Y$ and $T.$%
\end{tabular}
\end{equation*}
Let $M_{1,}M_{1}^{\prime },M_{2},M_{2}^{\prime }$ denote these four
monomials.

>From Remark \ref{remark3.11} it follows that all these monomials are
equipollent to each other and are equipollent to a certain monomial which we
will call $XYPZT$. Set 
\begin{eqnarray*}
Q_{1} &=&XYPZT,\qquad Q_{1}^{\prime }=XYPTZ, \\
Q_{2} &=&YXPTZ,\qquad Q_{2}^{\prime }=YXPZT.
\end{eqnarray*}

We are going to construct a series of simple chains which will constitute a
decomposition of $S^{\prime }$ as follows:\bigskip

\begin{tabular}{lll}
Name of the chain & Extreme positive monomial & Extreme negative monomial \\ 
\multicolumn{1}{c}{$C^{\prime \prime \prime }+C_{1\smallskip }$} & 
\multicolumn{1}{c}{$M_{1}=X\square T$} & \multicolumn{1}{c}{$-M_{1}^{\prime
}=-X\square Z$} \\ 
\multicolumn{1}{c}{$C_{1}^{\prime }\smallskip $} & \multicolumn{1}{c}{$%
M_{1}^{\prime }$} & \multicolumn{1}{c}{$-M_{2}=-Y\square Z$} \\ 
\multicolumn{1}{c}{$C_{2}\smallskip $} & \multicolumn{1}{c}{$M_{2}$} & 
\multicolumn{1}{c}{$-M_{2}^{\prime }=-Y\square T$} \\ 
\multicolumn{1}{c}{$C_{2}^{\prime }+C^{\prime \prime }\smallskip $} & 
\multicolumn{1}{c}{$M_{2}^{\prime }$} & \multicolumn{1}{c}{$-M_{1}$} \\ 
\multicolumn{1}{c}{$D_{1}\smallskip $} & \multicolumn{1}{c}{$M_{1}$} & 
\multicolumn{1}{c}{$-Q_{1}$} \\ 
\multicolumn{1}{c}{$D_{1}^{\prime }\smallskip $} & \multicolumn{1}{c}{$%
M_{1}^{\prime }$} & \multicolumn{1}{c}{$-Q_{1}^{\prime }$} \\ 
\multicolumn{1}{c}{$D_{2}\smallskip $} & \multicolumn{1}{c}{$M_{2}$} & 
\multicolumn{1}{c}{$-Q_{2}$} \\ 
\multicolumn{1}{c}{$D_{2}^{\prime }\smallskip $} & \multicolumn{1}{c}{$%
M_{2}^{\prime }$} & \multicolumn{1}{c}{$-Q_{2}^{\prime }$} \\ 
\multicolumn{1}{c}{$E_{1}\smallskip $} & \multicolumn{1}{c}{$Q_{1}$} & 
\multicolumn{1}{c}{$-Q_{1}^{\prime }$} \\ 
\multicolumn{1}{c}{$E_{1}^{\prime }\smallskip $} & \multicolumn{1}{c}{$%
Q_{1}^{\prime }$} & \multicolumn{1}{c}{$-Q_{2}$} \\ 
\multicolumn{1}{c}{$E_{2}\smallskip $} & \multicolumn{1}{c}{$Q_{2}$} & 
\multicolumn{1}{c}{$-Q_{2}^{\prime }$} \\ 
\multicolumn{1}{c}{$E_{2}^{\prime }$} & \multicolumn{1}{c}{$Q_{2}^{\prime }$}
& \multicolumn{1}{c}{$-Q_{1}^{{}}$}
\end{tabular}
\bigskip

\noindent We can suppose that every monomial of the chain $D_{1}$ has as
first factor $X$ and last factor $T$; thus $D_{1}$ is both a simple chain of
the first and second kind, and similarly for other $D$ and $D^{\prime }$
chains. Furthermore, we can suppose that the $E$ and $E^{\prime }$ chains
are reduced to a single trinomial product, for example, 
\begin{equation*}
E_{1}=XYP(ZT-TZ-[Z,T]).
\end{equation*}

The closed chain $S^{\prime }=(C^{\prime \prime \prime
}+C_{1})+C_{1}^{\prime }+C_{2}+(C_{2}^{\prime }+C^{^{\prime \prime }})$ can
be decomposed into five closed chains as follows: 
\begin{align*}
U_{1}& =\underset{C^{\prime \prime \prime }+C_{1}}{\underbrace{(M_{1}\cdots
-M_{1}^{\prime })}}+\underset{D_{1}^{\prime }}{\underbrace{M_{1}^{\prime
}\cdots -Q_{1}^{\prime }}}+\underset{-E_{1}}{\underbrace{XYP(TZ-ZT-[T,Z])}}+%
\underset{-D_{1}}{\underbrace{Q_{1}\cdots -M_{1}},} \\
U_{1}^{\prime }& =\underset{C_{1}^{\prime }}{\underbrace{M_{1}\cdots -M_{2}}}%
+\underset{D_{2}}{\underbrace{M_{2}\cdots -Q_{2}}}+\underset{-E_{1}^{\prime }%
}{\underbrace{(YX-XY-[Y,X])PTZ}}+\underset{-D_{1}^{\prime }}{\underbrace{%
Q_{1}^{\prime }\cdots -M_{1}^{\prime }}}, \\
U_{2}& =\underset{C_{2}}{\underbrace{M_{2}\cdots -M_{2}^{\prime }}}+%
\underset{D_{2}^{\prime }}{\underbrace{M_{2}^{\prime }\cdots -Q_{2}^{\prime }%
}}+\underset{-E_{2}}{\underbrace{YXP(ZT-TZ-[Z,T])}}+\underset{-D_{2}}{%
\underbrace{Q_{2}\cdots -M_{2}}}, \\
U_{2}^{\prime }& =\underset{C_{2}+C^{\prime \prime }}{\underbrace{%
M_{2}^{\prime }\cdots -M_{1}}}+\underset{D_{1}}{\underbrace{M_{1}\cdots
-Q_{1}}}+\underset{-E_{2}^{\prime }}{\underbrace{(XY-YX-[X,Y])PZT}}+%
\underset{-D_{2}^{\prime }}{\underbrace{Q_{2}^{\prime }\cdots -M_{2}^{\prime
}}}, \\
V& =\underset{E_{1}}{\underbrace{XYP(ZT-TZ-[Z,T])}}+\underset{E_{1}^{\prime }%
}{\underbrace{(XY-YX-[X,Y])PTZ}} \\
& \qquad \qquad +\underset{E_{2}}{\underbrace{(YXP(TZ-ZT-[T,Z]})}+\underset{%
E_{2}^{\prime }}{\underbrace{(YX-XY-[Y,X])PZT}}.
\end{align*}
Clearly, $U_{1}+U_{1}^{\prime }+U_{2}+U_{2}^{\prime }+V=(C^{\prime \prime
\prime }+C_{1})+C_{1}^{\prime }+C_{2}+(C_{2}^{\prime }+C^{\prime \prime
})=S. $ We must show that each one of the five closed chains above is the
head of an identically zero regular sum. The first four chains are of the
form 
\begin{equation*}
U_{1}=XH_{1},\quad U_{1}^{\prime }=H_{1}^{\prime }Z,\quad U_{2}=YH_{2},\quad
U_{2}^{\prime }=H_{2}^{\prime }T,
\end{equation*}
where each chain $H_{1},H_{1}^{\prime },H_{2},H_{2}^{\prime }$ is a closed
chain of degree $p-1$; therefore by induction, each is the head of an
identically zero regular sum. It follows that $U_{1}$, $U_{1}^{\prime }$, $%
U_{2}$, and $U_{2}^{\prime }$ are identically zero, and therefore each of
them can be considered as the head of an identically zero regular sum of
degree $p$.

Finally for $V$, it is the head of the sums 
\begin{multline*}
XYP(ZT-TZ-[Z,T])+(XY-YX-[X,Y])PTZ \\
-YXP(-TZ+ZT-[Z,T])-(-YX+XY-[X,Y])PZT \\
-[X,Y]P(ZT-TZ-[Z,T])-(XY-YX-[X,Y])P[T,Z],
\end{multline*}
which can be expanded and rearranged as 
\begin{multline*}
XYPZT-XYPTZ+XYPTZ-YXPTZ+YXPTZ \\
-YXPZT+YXPZT-XYPZT-XYP[Z,T]-[X,Y]PTZ \\
+YXP[Z,T]+[X,Y]PZT-[X,Y]PZT+[X,Y]PTZ \\
-XYP[T,Z]+YXP[T,Z]+[X,Y]P[Z,T]+[X,Y]P[T,Z],
\end{multline*}
which is identically zero. Since $0$ is a regular sum, it follows that $V$
is the head of an identically zero regular sum of degree $p$.

Note that our analysis remains unchanged when two or more of the operators $%
X $, $Y$, $Z$, $T$ are identical. For example, when $X=Y$, then $%
E_{1}^{\prime }=E_{2}^{\prime }=0$, and we set $Q_{1}=Q_{2}^{\prime
}=X(XP)ZT=XP^{\prime }ZT$, $Q_{2}=Q_{1}^{\prime }=X(XP)TZ=XP^{\prime }TZ$.
The definition of the various chains remains the same, and we can
immediately verify that $V$ is identically zero. Finally, in order that this
proof is valid, $p$ must be greater than three since the chain $V$ must have
at least four factors. But this was the assumption in our inductive
hypothesis. Thus the proof of Lemma \ref{lemma3.9}, and hence of Theorem \ref
{theorem3.8}, is achieved.%
%TCIMACRO{
%\TeXButton{End Proof}{\endproof%
%}}%
%BeginExpansion
\endproof%
%
%EndExpansion

\begin{corollary}
\label{corollary3.12}$($The so-called Birkhoff-Witt Theorem.$)$ \textit{Let }%
$\mathcal{U}(\mathcal{L})$\textit{\ denote the universal enveloping algebra
of a Lie algebra over a }$($\textit{commutative}$)$\textit{\ field of
characteristic zero. If }$\{X_{1},\ldots ,X_{n}\}$\textit{\ is a basis of }$%
\mathcal{L}$\textit{\ and if }$(\alpha )=(\alpha _{1},\ldots ,\alpha _{n})$%
\textit{\ denotes an }$n$\textit{-tuple of integers }$\geq 0$\textit{, set }$%
X^{(\alpha )}=X_{1}^{\alpha _{1}}\cdots X_{n}^{\alpha _{n}},\tilde{X}%
^{(\alpha )}=\limfunc{Sym}(X^{\alpha }),\left| (\alpha )\right| =\alpha
_{1}+\cdots +\alpha _{n}$. \textit{Then the set} $\{\tilde{X}^{(\alpha
)}\}_{(\alpha )}$, \textit{for all }$(\alpha )$\textit{\ such that }$\left|
(\alpha )\right| \geq 0$, \textit{forms a vector space basis for }$\mathcal{U%
}(\mathcal{L)}$\textit{. Moreover, any set of elements of }$\mathcal{U}(%
\mathcal{L)}$\textit{\ of the form }$\{X_{i_{1}}\cdots X_{i_{m}},1\leq
i_{j}\leq n,1\leq j\leq m,m\geq 0\}$\textit{, where each }$X_{i_{1}}\cdots
X_{i_{m}}$\textit{\ is a representative of an equipollence class }$%
X^{(\alpha )}$\textit{\ for all distinct }$(\alpha )$\textit{\ such that }$%
\left| (\alpha )\right| \geq 0$\textit{, is a basis of }$\mathcal{U}(%
\mathcal{L})$\textit{; in particular, the set of ordered monomials }$%
\{X^{(\alpha )},\left| (\alpha )\right| \geq 0\}$\textit{\ forms a basis for 
}$\mathcal{U}(\mathcal{L})$.
\end{corollary}

\begin{proof}
>From our discussion pertaining to Poincar\'{e}'s discovery of the universal
enveloping algebra of a Lie algebra, it follows that the quotient algebra of
the polynomial algebra $\mathcal{A}$ modulo the equivalence relation $\sim $
can be regarded as the universal enveloping algebra of the Lie algebra $%
\mathcal{L}$ generated by $X_{1},\ldots ,X_{n}$.

Define a map from $\mathcal{A}$ to $\mathcal{R}$, the vector space of all
regular polynomials in $\mathcal{A}$, by assigning to each polynomial $A$ in 
$\mathcal{A}$ the unique regular polynomial $\bar{A}$ equivalent to $A$ as
defined by Theorem \ref{theorem3.8}. From the proof of Theorem \ref
{theorem3.8} it follows that this map is linear, and that it is surjective
since the unique regular polynomial equivalent to a given regular polynomial
is itself. The kernel of this homomorphism is, by the definition of the
equivalence relation $\sim $, the vector space spanned by all trinomials of
the form $P(XY-YX-[X,Y])Q$ for arbitrary $P$ and $Q$ in $\mathcal{A}$. Let $%
\mathcal{I}$ denote this kernel, then obviously $\mathcal{I}$ is a two-sided
ideal of $\mathcal{A}$. It follows from the first isomorphism theorem that $%
\mathcal{A}/\mathcal{I}$ is isomorphic to $\mathcal{R}$ as vector spaces.
>From the remark following Lemma \ref{lemma3.5}, it follows that the set $\{%
\tilde{X}^{(\alpha )}\}_{(\alpha )},$ $\left| (\alpha )\right| \geq 0$,
forms a basis for $\mathcal{R}$, and hence a basis for $\mathcal{U}(\mathcal{%
L})\cong \mathcal{A}/\mathcal{I}$ via the isomorphism above. Note that we
have shown following Lemma \ref{lemma3.5} that $\mathcal{P}$ is isomorphic
to $\mathcal{R}$ via the isomorphism $\Phi $, therefore $\mathcal{P}$ is
isomorphic to $\mathcal{U}(\mathcal{L})$. For the second part of the
theorem, we remark that it follows from Theorem \ref{Theorem3.3} that each $%
X_{i_{1}}\cdots X_{i_{m}}$ is equipollent to a unique regular polynomial $%
\tilde{X}^{(\alpha )}$ for some $(\alpha )\in \Lambda _{m}$. Thus it
suffices to consider the set $\{X^{(\alpha )},\left| (\alpha )\right| \geq
0\}$. We also remark that it suffices to show that the set $\{X^{(\alpha
)},(\alpha )\in \Lambda _{m}\}$ is linearly independent in $\mathcal{U}(%
\mathcal{L})$ for all $m\geq 0$, since $\mathcal{U}(\mathcal{L})$ is a
filtered algebra. From the proof of Theorem \ref{theorem3.8} it follows that
each $X^{(\alpha )},(\alpha )\in \Lambda _{m}$, is equivalent to a unique
regular polynomial of the form $\tilde{X}^{(\alpha )}+P_{(\alpha )}$, where $%
P_{(\alpha )}$ is a regular polynomial of degree $<m$. Thus if for some
scalars $\lambda _{(\alpha )}\in \mathbb{K}$ such that $\sum_{(\alpha )\in
\Lambda _{m}}\lambda _{(\alpha )}X^{(\alpha )}$ is zero in $\mathcal{U}(%
\mathcal{L})$ (i.e., equivalent to $0$), then since the equivalence relation 
$\sim $ is linear it follows that the regular polynomial $\sum_{(\alpha )\in
\Lambda _{m}}\lambda _{\alpha }(\tilde{X}^{(\alpha )}+P_{\alpha })$ is
equivalent to $0$. It follows from Theorem \ref{theorem3.8} (or more
precisely Lemma \ref{lemma3.9}) that 
\begin{equation*}
\sum_{(\alpha )\in \Lambda _{m}}\lambda _{(\alpha )}\tilde{X}^{(\alpha
)}+\sum_{(\alpha )\in \Lambda _{m}}\lambda _{(\alpha )}P_{(\alpha )}
\end{equation*}
must be identically zero. Since $d^{0}(\tilde{X}^{(\alpha )})=m$ and $%
d^{0}(P_{(\alpha )})<m$ for all $(\alpha )\in \Lambda _{m}$, it follows that 
\begin{equation*}
\sum_{(\alpha )\in \Lambda _{m}}\lambda _{(\alpha )}\tilde{X}^{(\alpha )}=0;
\end{equation*}
and hence by the first part of the proof of the theorem, it follows that $%
\lambda _{(\alpha )}=0$ for all $(\alpha )\in \Lambda _{m}$. This completes
the proof of the theorem.
\end{proof}

\begin{remark}
\textrm{\label{remark3.13}}We note that throughout this section the basis
for the Lie algebra $\mathcal{L}$ can be the \textit{infinite} set $%
\{X_{1},\ldots ,X_{n},\ldots \}$. The tensor algebra $\mathcal{T}$ remains
isomorphic to the non-commutative algebra $\mathcal{A}$ of polynomials in
infinitely many variables $X_{1},\ldots ,X_{n},\ldots $ (see, e.g., 
%TCIMACRO{
%\TeXButton{TeX field}{\cite[Prop.~(2.4), p.~ 40]{Schwartz}%
%}}%
%BeginExpansion
\cite[Prop.~(2.4), p.~ 40]{Schwartz}%
%
%EndExpansion
), and every argument remains the same. As a special case, let $\mathcal{L}%
_{n}$ (resp.$\ \mathcal{A}_{n}$) denote the Lie algebra (resp.\ the
non-commutative polynomial algebra) generated by $X_{1},\ldots ,X_{n}$. Let $%
\mathcal{L}$ (resp.\ $\mathcal{A}$) denote the \textit{inductive limit} of $%
\mathcal{L}_{n}$ (resp.$\ \mathcal{A}_{n}$); then all theorems in this
section can be easily generalized. For example, let 
\begin{align*}
& X_{ij}=x_{i}\frac{\partial }{\partial x_{j}},1\leq i,j\leq n;\text{ then}
\\
& [X_{ij},X_{kl}]=\delta _{jk}X_{il}-\delta _{li}X_{kj},
\end{align*}
and $\{X_{ij},1\leq i,j\leq n\}$ generates the Lie algebra $\frak{g}\frak{l}%
_{n}$ and the associative polynomial algebra $\mathcal{A}_{n}$,
respectively. By letting $n\rightarrow \infty $ we get the Lie algebra $%
\frak{g}\frak{l}_{\infty }$ and $\mathcal{A}_{\infty }$, respectively.
Another example is the Heisenberg Lie algebra $\mathcal{H}_{n}$ spanned by
the vector fields 
\begin{equation*}
P_{j}=\frac{i\sqrt{2}}{2}\left( x_{j}+\frac{\partial }{\partial x_{j}}%
\right) ,\quad Q_{j}=\frac{\sqrt{2}}{2}\left( -x_{j}+\frac{\partial }{%
\partial x_{j}}\right) ,\quad 1\leq j\leq n,\quad i=\sqrt{-1},
\end{equation*}
and $R=iI$, where $I$ is the identity operator. Then we have the commutation
relations 
\begin{equation*}
\lbrack
P_{j},P_{k}]=[Q_{j},Q_{k}]=0,\;[P_{j},R]=[Q_{j},R]=0,\;[P_{j},Q_{k}]=-\delta
_{jk}R,\quad 1\leq j,k\leq n.
\end{equation*}
Let $\mathcal{A}_{n}$ denote the algebra of non-commutative polynomials in
the vector fields $P_{1},\ldots ,P_{n}$, $Q_{1},\ldots ,Q_{n}$ and $R$. Then
when $n\rightarrow \infty $ we obviously have the generalization of the
theorems in this section to the Heisenberg Lie algebra $\mathcal{H}_{\infty
} $, and hence to $\mathcal{A}_{\infty }$.
\end{remark}

\section{Conclusion\label{Conclusion}}

A. Einstein said ``A good idea is very rare''. We reckon that there are, at
the very least, three ``good'' ideas in \cite{Po1}, namely, the universal
enveloping algebra of a Lie algebra, the symmetrization map, and the proof
of the so-called Birkhoff-Witt theorem. And in our opinion, none of these
were properly appreciated and recognized. We have gone to great length and
sometimes with repetitive arguments to try to convince the mathematics
community of what a great feat Poincar\'{e} has achieved in \cite{Po1}. But
even if we fail, we would be much wiser by our reading a masterpiece by a
great master.

We leave the reader with the following thought of Paul Painlev\'{e}, another
great master, in the obituary written for the newspaper \emph{Le Temps} (and
reprinted in \cite{Pa}), on July 18, 1912, the day after Poincar\'{e} died:

``Henri Poincar\'{e} n'a pas \'{e}t\'{e} seulement un grand cr\'{e}ateur
dans les sciences positives. Il a \'{e}t\'{e} un grand philosophe et un
grand \'{e}crivain. Certains de ses aphorismes font songer \`{a} Pascal: `La
pens\'{e}e n'est qu'un \'{e}clair entre deux longues nuits, mais c'est cet
\'{e}clair qui est tout'. Son style traduit la d\'{e}marche m\^{e}me de sa
pens\'{e}e: des formules br\`{e}ves et saisissantes, parodoxales parfois
quand on les isole, r\'{e}unies par des explications h\^{a}tives, qui
rejettent des d\'{e}tails faciles pour ne dire que l'essentiel. C'est
pourquoi des critiques superficiels lui ont reproch\'{e} d'\^{e}tre
`d\'{e}cousu': la v\'{e}rit\'{e}, c'est que, sans \'{e}ducation scientifique
pr\'{e}alable, une telle d\'{e}marche logique est difficile \`{a}
\'{e}galer: le lion ne fait pas des enjamb\'{e}es de souris''\footnote{%
This can be roughly translated as follows: ``Henri Poincar\'{e} was not only
a great creator in the positive sciences. He was a great philosopher and a
great writer. Some of his aphorisms make us think of Pascal: `Thought is
just a flash of lightning in the middle of two long nights, but it is this
lightning that is everything'. His style reflects the very development of
his thought: brief and startling formulae, sometimes paradoxical when one
isolates them, joined together by hasty explanations, which reject easy
details in order just to express the essential. That is why superficial
critiques reproach him as being `incoherent': the truth is that, without
prerequisite education, such logical development is difficult to match: the
lion does not take a mouse's paces.''}.\medskip

\begin{acknowledgement}
The authors wish to thank the editors of the \emph{Revue d'Histoire des
Math\'{e}matiques} whose comments helped to improve the overall quality of
this article. They also wish to thank the technical staff of the Department
of Mathematics of the University of Iowa, Mr. Brian Treadway and Ms. Cymie
Wehr, for their impeccable job and for their patience in preparing this
manuscript.
\end{acknowledgement}

\newpage 

\section*{References}

\begin{list}{}{\setlength{\leftmargin}{\customleftmargin}
\setlength{\labelwidth}{\customleftmargin}
\addtolength{\labelwidth}{-\labelsep}
\setlength{\itemsep}{0.5\customskipamount}
\setlength{\parsep}{0.5\customskipamount}
\def\makelabel#1{\upshape#1\hss}
\frenchspacing}

\item[\textit{Encyclopaedia}]
\item[\kern\customkern\lbrack 1988--1994\rbrack ]  \textit{Encyclopaedia of
Mathematics},
Editor-in-Chief I.M. Vinogradov,
Dordrecht-Boston:\ Kluwer\ Academic Publishers, 1988--1994,
translated from {\cyr Matematicheskaya \`Entsiklopediya,
glavny\u\i{} redaktor I. M. Vinogradov,
Tom 1, 495, Moskva: Sovet\cydot skaya \`Entsiklopediya, 1977}.

\item[\textsc{Barrow-Green (J.)}]
\item[\kern\customkern\lbrack 1997\rbrack ]    \textit{Poincar\'{e} and
the Three-Body Problem}, Providence:\ American Mathematical Society, 1997.

\item[\textsc{Bell (E.T.)}]
\item[\kern\customkern\lbrack 1937\rbrack ]  \textit{Men of Mathematics}, New York:
Simon \& Schuster, 1937.

\item[\textsc{Birkhoff (Garrett)}]
\item[\kern\customkern\lbrack 1937\rbrack ]  Representability of Lie
algebras and Lie groups by matrices, \textit{Ann.\ of Math.}, 38 (April
1937), pp.\ 526--532.

\item[\textsc{Bourbaki (N.)}]
\item[\kern\customkern\lbrack 1960\rbrack ]  \textit{\'{E}l\'{e}ments de
Math\'ematiques, Groupes et Alg\'{e}bres de Lie}, Chap.\ 1, Paris:\ Hermann,
1960.

\item[\kern\customkern\lbrack 1969\rbrack ]  \textit{\'{E}l\'{e}ments
d'Histoire des
Math\'ematiques}, Paris:\ Hermann, 1969.

\item[\kern\customkern\lbrack 1972\rbrack ]  \textit{\'{E}l\'{e}ments de
Math\'ematiques, Fasc.\ XXXVII, Groupes et Alg\`{e}bres de Lie}, Paris:\
Hermann, 1972.

\item[\kern\customkern\lbrack 1975\rbrack ]  \textit{Elements of Mathematics,
Lie Groups and Lie Algebras, }Part I, Chapter I-3\textit{, }Reading, MA:
Addison-Wesley, 1975; Paris:\ Hermann,\textit{\ }1971--73.

\item[\textsc{Boyer (C.)}]
\item[\kern\customkern\lbrack 1968\rbrack ]  \textit{A History of Mathematics},
Princeton:\ Princeton University Press, 1968.

\item[\textsc{Cartan (H.) \& Eilenberg (S.)}]
\item[\kern\customkern\lbrack 1956\rbrack ]  
\textit{Homological Algebra}, Princeton:\ Princeton University Press, 1956.

\item[\textsc{Chevalley (C.)}]
\item[\kern\customkern\lbrack 1955\rbrack ]  \textit{Th\'{e}orie des Groupes
de Lie}, vol.\ III, Paris:\ Hermann, 1955.

\item[\textsc{Cohn (P.M.)}]
\item[\kern\customkern\lbrack 1981\rbrack ]  \textit{Universal Algebra},\emph{\ }%
Dordrecht:\ Reidel, 1981.

\item[\textsc{Dixmier (J.)}]
\item[\kern\customkern\lbrack 1974\rbrack ]  \textit{Alg\`{e}bres
Enveloppantes}%
, Paris:\ Gauthier-Villars, 1974.

\item[\textsc{Gittleman (A.)}]
\item[\kern\customkern\lbrack 1975\rbrack ]  \textit{History of Mathematics}%
, Columbus, OH:\ C.E. Merrill Publishing Company, 1975.

\item[\textsc{Godement (R.)}]
\item[\kern\customkern\lbrack 1982\rbrack ]  \textit{Introduction \`a la
Th\'{e}orie des Groupes de Lie}, 2 vols., Publ.\ Math.\ Univ.\ Paris VII,
Paris:\ Univ.\ Paris VII, 1982.

\item[\textsc{Harish-Chandra}]
\item[\kern\customkern\lbrack 1949\rbrack ]  On representations of
Lie algebras, \textit{Ann.\ of Math.}, 50 (October 1949), pp.\ 900--915.

\item[\kern\customkern\lbrack 1951\rbrack ]  On some applications of
the universal enveloping algebra of a semisimple Lie algebra, \textit{%
Trans.\ Amer.\ Math.\ Soc.}, 70--71 (1951), pp.\ 28--96.

\item[\textsc{Hoffman (K.) \& Kunze (R.)}]
\item[\kern\customkern\lbrack 1971\rbrack ]  \textit{Linear algebra}%
, second ed., Englewood Cliffs, NJ:\ Prentice-Hall, 1971.

\item[\textsc{Humphreys (J.E.)}]
\item[\kern\customkern\lbrack 1972\rbrack ]  \textit{Introduction to Lie
Algebras and Representation Theory}, New York:\ Springer-Verlag, 1972.

\item[\textsc{Jacobson (N.)}]
\item[\kern\customkern\lbrack 1962\rbrack ]  \textit{Lie Algebras}, New York:
John Wiley \& Sons, 1962.

\item[\textsc{Knapp (A.)}]
\item[\kern\customkern\lbrack 1986\rbrack ]  \textit{Representation Theory of
Semisimple Group}, Princeton:\ Princeton University Press, 1986.

\item[\textsc{Kourotchkine}]
\item[\kern\customkern\lbrack 1951\rbrack ]  \textit{Mat.\ Sb.}, 28, (70),
2, (1951), pp.\ 467--472.

\item[\textsc{Kuro\'{s} (A.G.)}]
\item[\kern\customkern\lbrack 1963\rbrack ]  \textit{Lectures on General
Algebras}, New York:\ Chelsea, 1963.

\item[\textsc{Lang (S.)}]
\item[\kern\customkern\lbrack 1965\rbrack ]  \textit{Algebra}, Reading, MA:
Addison-Wesley, 1965.

\item[\textsc{Lazard (M.)}]
\item[\kern\customkern\lbrack 1952\rbrack ]  Sur les algebres enveloppantes
universelles de certaines alg\`{e}bres de Lie, \textit{C.R. Acad.\ Sci.\
Paris S\'{e}r.\ I} (18 Feb.\ 1952).

\item[\kern\customkern\lbrack 1954\rbrack ]  Sur les alg\`{e}bres enveloppantes
de certaines alg\`{e}bres de Lie, \textit{Publ.\ Sci.\ Univ.\ Alger., Ser.\ A},
1 (1954), pp.\ 281--294.

\item[\textsc{Lorentz (H.A.)}]
\item[\kern\customkern\lbrack 1921\rbrack ]  Deux m\'{e}moires de Henri
Poincar\'{e} sur la physique math\'{e}matique, \textit{Acta Math.}, 38 (1921),
pp.\ 293--308.

\item[\textsc{Painlev\'{e} (P.)}]
\item[\kern\customkern\lbrack 1921\rbrack ]  Henri Poincar\'{e}, \textit{%
Acta Math.}, 38 (1921), pp.\ 309--402.

\pagebreak

\item[\textsc{Poincar\'{e} (H.)}]
\item[\kern\customkern\lbrack 1881\rbrack ]  Formes cubiques
ternaires et quaternaires, \textit{Jour. de l'\'{E}cole Polytechnique
Paris}, XXXI (1881), pp.\ 199--253.

\item[\kern\customkern\lbrack 1883\rbrack ]  Sur la reproduction des
formes,\textit{\ Comptes Rendus des S\'{e}ances de l'Acad.\ des\
Sci.\ Paris} (29 Oct.\ 1883)\ 

\item[\kern\customkern\lbrack 1899\rbrack ]  Sur les groupes
continus, \textit{C.R.\ Acad.\ Sci.,} 128 (1899), pp.\ 1065--1069.

\item[\kern\customkern\lbrack 1900\rbrack ]  Sur les groupes
continus, \textit{Trans.\ Cambr.\ Philos.\ Soc.,} 18 (1900), pp.\
220--255 = \textit{Oeuvres de Henri Poincar\'{e}}, vol., III, Paris:
Gauthier-Villars, (1934), pp.\ 173--212 (1899).

\item[\kern\customkern\lbrack 1901\rbrack ]  Quelques remarques sur
les groupes continus, \textit{Rend. Circ. Mat. Palermo},
15 (1901).

\item[\kern\customkern\lbrack 1906\rbrack ]  Sur la dynamique de
l'election, \textit{Circolo Matematico di Palermo}, 21 (1906).

\item[\kern\customkern\lbrack 1908\rbrack ]  Nouvelles remarques sur
les groupes continus, \textit{Rend. Circ. Mat. Palermo},
25 (1908)

\item[\kern\customkern\lbrack 1912\rbrack ]  Sur la th\'{e}orie des
quanta, \textit{Journal de Physique Th\'{e}orique et Appliqu\'{e}e
5e S\'{e}r.}, 2 (1912).

\item[\kern\customkern\lbrack 1916--\rbrack ]  \textit{Oeuvres de Henri
Poincar\'{e}}, 11 vols., Paris:\ Gauthier-Villars, \linebreak\mbox{1916--\quad .}

\item[\textit{In memoriam}]
\item[\kern\customkern\lbrack 1921\rbrack ]  Henri Poincar\'{e} in memoriam, \textit{%
Acta Math.}, 38 (1921).

\item[\textsc{Schmid (W.)}]
\item[\kern\customkern\lbrack 1982\rbrack ]  Poincar\'{e} and Lie
groups, \textit{Bull. Amer. Math. Soc. (N.S.)},
6 (1982), pp.\ 175--186.
This article is also reprinted in
\textit{The Mathematical Heritage of Henri Poincar\'{e}},
Proceedings of Symposia in Pure Mathematics, Vol.~39,
edited by Felix E. Browder, Providence, R.I.:\
American Mathematical Society, 1983.

\item[\textsc{Schwartz (L.)}]
\item[\kern\customkern\lbrack 1998 (1975)\rbrack ]  \textit{Les
tenseurs}, Paris:\ Hermann,
1998; 1e \'{e}d.\ 1975.

\item[\textsc{Varadarajan (V.S.)}]
\item[\kern\customkern\lbrack 1984 (1974)\rbrack ]  \textit{Lie
Groups, Lie Algebras and Their Representations}, New York:\ Springer-Verlag,
1984; Englewood Cliffs, NJ:\ Prentice-Hall, 1974.

\item[\textsc{Weyl (H.)}]
\item[\kern\customkern\lbrack 1946\rbrack ]  \textit{The Classical Groups, Their
Invariants and Representations}, second ed., Princeton:\ Princeton Univ.\
Press, 1946.

\item[\textsc{Witt (E.)}]
\item[\kern\customkern\lbrack 1937\rbrack ]  Treue Darstellung Liescher Ringe, \textit{%
J. Reine Angew. Math.}, 177 (1937), pp.\
152--160.
\end{list}
\pagebreak \textbf{PLEASE DISREGARD THE FOLLOWING PAGES}


\begin{thebibliography}{Poincar\'{e} in memoriam 1921}
\bibitem[Barrow-Green 1997]{B-G}  Barrow-Green, J., \textit{Poincar\'{e} and
the Three-Body Problem}, Providence: American Mathematical Society, 1997.

\bibitem[Bell 1937]{Be}  Bell, E.T., \textit{Men of Mathematics}, New York:
Simon \& Schuster, 1937.

\bibitem[Birkhoff 1937]{Bi}  Birkhoff, Garrett, Representability of Lie
algebras and Lie groups by matrices, \textit{Ann. of Math.}, 38 (April
1937), pp. 526--532.

\bibitem[Bourbaki 1972]{Bo1}  Bourbaki N., \textit{\'{E}l\'{e}ments de
mathematics, Fasc. XXXVII, Groupes et Alg\`{e}bras de Lie}, Paris: Hermann,
1972.

\bibitem[Bourbaki 1975]{Bo2}  Bourbaki, N., \textit{Elements of Mathematics,
Lie Groups and Lie Algebras, }Part I, Chapter I-3\textit{, }Reading, MA:
Addison-Wesley, 1975; Paris: Hermann,\textit{\ }1971-73.

\bibitem[Bourbaki 1960]{Bo3}  Bourbaki, N., \textit{\'{E}l\'{e}ments de
Mathematics, Groupes et alg\'{e}bra de Lie}, Chap. 1, Paris: Hermann, 1960.

\bibitem[Bourbaki 1969]{Bo}  Bourbaki, N., \'{E}l\'{e}ments d'histoire des
math\'{e}matiques, Paris: Hermann, 1969.

\bibitem[Boyer 1968]{By}  Boyer, C., \textit{A History of Mathematics},
Princeton: Princeton University Press, 1968.

\bibitem[Cartan \& Eilenberg 1956]{C-E}  Cartan, H. \& Eilenberg, S., 
\textit{Homological Algebra}, Princeton: Princeton University Press, 1956.

\bibitem[Chevalley 1955]{Ch}  Chevalley, C., \textit{Th\'{e}orie des groupes
de Lie}, vol. III, Paris: Hermann, 1955.

\bibitem[Cohn 1981]{Co}  Cohn, P.M., \textit{Universal algebra},\emph{\ }%
Dordrecht: Reidel, 1981.

\bibitem[Dixmier 1974]{Di}  Dixmier, J., \textit{Alg\`{e}bres enveloppantes}%
, Paris: Gauthier-Villars, 1974.

\bibitem[Encyclopaedia 1988-1994]{En}  \textit{Encyclopaedia of Mathematics}%
, Dordrecht-Boston: Kluwer\ Academic Publishers, 1988-1994.

\bibitem[Gittleman 1975]{Gi}  Gittleman, A., \textit{History of Mathematics}%
, Columbus, OH: C.E. Merrill Publishing Company, 1975.

\bibitem[Godement 1982]{Go}  Godement, R., \textit{Introduction a la
th\'{e}orie des groupes de Lie}, 2 vols., Publ.\ Math.\ Univ.\ Paris VII,
Paris: Univ. Paris VII, 1982.

\bibitem[Harish-Chandra 1949]{H-C1}  Harish-Chandra, On representations of
Lie algebras, \textit{Ann. of Math.}, 50 (October 1949), pp. 900--915.

\bibitem[Harish-Chandra 1951]{H-C2}  Harish-Chandra, On some applications of
the universal enveloping algebra of a semisimple Lie algebra, \textit{%
Trans.\ Amer.\ Math.\ Soc.}, 70--71 (1951), pp. 28--96.

\bibitem[Hoffman \& Kunze 1971]{H-K}  Hoffman, K. \& Kunze R., \textit{%
Linear algebra}, second ed., Englewood Cliffs, NJ: Prentice-Hall, 1971.

\bibitem[Humphreys 1972]{Hu}  Humphreys, J.E., \textit{Introduction to Lie
Algebras and Representation Theory}, New York: Springer-Verlag, 1972.

\bibitem[Jacobson 1962]{Ja}  Jacobson, N., \textit{Lie Algebras}, New York:
John Wiley \& Sons, 1962.

\bibitem[Knapp 1986]{Kn}  Knapp, A., \textit{Representation Theory of
Semisimple Group}, Princeton: Princeton University Press, 1986.

\bibitem[Kourotchkine 1951]{Kr}  Kourotchkine, \textit{Mat. Sb.}, 28, (70),
2, (1951), pp. 467--472.

\bibitem[Kuros 1963]{Ku}  Kuro\'{s}, A.G., \textit{Lectures on General
Algebras}, New York: Chelsea, 1963.

\bibitem[Lang 1965]{La}  Lang, S., \textit{Algebra}, Reading, MA:
Addison-Wesley, 1965.

\bibitem[Lazard 1952]{Lz-1}  Lazard, M., Sur les algebres enveloppantes
universelles de certaines alg\`{e}bres de Lie, \textit{C.R.\ Acad.\ Sci.\
Paris\ S\'{e}r I} (18 Feb. 1952).

\bibitem[Lazard 1954]{Lz-2}  Lazard, M., Sur les alg\`{e}bres enveloppantes
de certaines alg\`{e}bres de Lie, \textit{Publ.\ Sci.\ Univ.\ Alger, Ser.A},
1 (1954), pp. 281-294.

\bibitem[Lorentz 1921]{Lo}  Lorentz, H.A., Deux m\'{e}moires de Henri
Poincar\'{e} sur la physique math\'{e}matique, \textit{Acta Math.} 38 (1921)
pp.\ 293--308.

\bibitem[Painlev\'{e} 1921]{Pa}  Painlev\'{e}, P., Henri Poincar\'{e}, 
\textit{Acta Math.} 38 (1921), pp. 309-402.

\bibitem[Poincar\'{e} 1900]{Po1}  Poincar\'{e}, H., Sur les groupes
continus, \textit{Trans.\ Cambr.\ Philos.\ Soc.,} 18\textbf{\ }(1900), pp.\
220--255 = \textit{Oeuvres de Henri Poincar\'{e}}, vol., III, Paris:
Gauthier-Villars, (1934), 173--212 (1899).

\bibitem[Poincar\'{e} 1916--]{Po2}  Poincar\'{e}, H, \textit{Oeuvres de
Henri Poincar\'{e}}, 11 vols., Paris: Gauthier-Villars, 1916--.

\bibitem[Poincar\'{e} in memoriam 1921]{Po3}  Henri Poincar\'{e} in
memoriam, \textit{Acta Math.}, 38 (1921).

\bibitem[Poincar\'{e} 1899]{Po4}  Poincar\'{e}, H., Sur les groupes
continus, \textit{C.R.\ Acad.\ Sci.,} 128 (1899), pp.1065-1069.

\bibitem[Poincar\'{e} 1881]{Po5}  Poincar\'{e}, H., Formes cubiques
ternaires et quaternaires, \textit{Jour de \'{E}'cole Polytechnique} \textit{%
Paris}, XXXI (1881), pp.\ 199--253.

\bibitem[Poincar\'{e} 1883]{Po6}  Poincar\'{e}, H., Sur la reproduction des
formes,\textit{\ Comptes rendus }des s\'{e}ances de l'\textit{Acad.\ des\
Sci.}\ \textit{Paris}, (29 Oct. 1883)\ 

\bibitem[Poincar\'{e} 1906]{Po7}  Poincar\'{e}, H., Sur la dynamique de
l'election, \textit{Circolo Matematico di Palermo}, 21 (1906).

\bibitem[Poincar\'{e} 1912]{Po8}  Poincar\'{e}, H., Sur la th\'{e}orie des
quarta, \textit{Journal de Physique Th\'{e}orique et Applique\'{e}} \textit{%
5\'{e}} \textit{Ser.}, 2 (1912).

\bibitem[Poincar\'{e} 1901]{Po9}  Poincar\'{e}, H., Quelques remarques sur
les groupes continus, \textit{Rendiconti del Circolo Matematico di Palermo},
15 (1901).

\bibitem[Poincar\'{e} 1908]{Po10}  Poincar\'{e}, H., Nouvelles remarques sur
les groupes continus, \textit{Rendiconti del Circolo Matematico di Palermo},
25 (1908)\emph{\ }

\bibitem[Varadarajan 1984 (1974)]{Va}  Varadarajan, V.S., \textit{Lie
groups, Lie algebras and their representations}, New York: Springer-Verlag,
1984; Englewood Cliffs, NJ: Prentice-Hall, 1974.

\bibitem[Schmid 1982]{W}  Poincar\'{e} and Lie Groups, Bulletin Amer. Math.
Soc. (NS) 6 (1982), pp. 175--186. This article is also reprinted in \emph{%
The Mathematical Heritage of Henri Poincar\'{e}}, Proceedings of \_\_\_\_ in
Pure Mathematics, Vol 39 (1983), edited by Felix Browder, Amer. Math. Soc.,
Providence, RI, USA.

\bibitem[Schwartz 1998 (1975)]{Schwartz}  Schwartz

\bibitem[Weyl 1946]{We}  Weyl, H., \textit{The classical groups, their
invariants and representations}, second ed., Princeton: Princeton Univ.\
Press, 1946.

\bibitem[Witt 1937]{Wi}  Witt, E., Treue Darstellung Liescher Ringe, \textit{%
Jour. f\"{u}r die reine und angewandte Mathematik}, 177 (1937), pp.\
152--160.
\end{thebibliography}
\end{document}